\documentclass[a4paper,12pt]{article}
\usepackage{amsmath}
\usepackage{amsfonts}
\usepackage{amssymb}
\usepackage{latexsym}
\usepackage{epsfig}
\usepackage{graphicx}
\usepackage{oldgerm}
\usepackage{theorem}

\def\bE{{\bf E}}

\def\sG{{\cal G}}

\def\1{{\bf 1}}

\def\sc{{\sf c}}
\def\cC{{\cal C}}
\def\cS{{\cal S}}

\def\mM{\mathfrak{M}}

\def\Det{{\rm Det}}

\theorembodyfont{\itshape}

\newtheorem{theorem}{Theorem}[section]
\newtheorem{lemma}[theorem]{Lemma}
\newtheorem{corollary}[theorem]{Corollary}
\newtheorem{proposition}[theorem]{Proposition}

\newcommand{\mib}[1]{\mbox{\boldmath $#1$}}
\newcommand{\SSC}[1]{\section{#1}\setcounter{equation}{0}}

\begin{document}

\title{\bf Complex Brownian Motion Representation \\
of the Dyson Model}
\author{
Makoto Katori
\footnote{
Department of Physics,
Faculty of Science and Engineering,
Chuo University, 
Kasuga, Bunkyo-ku, Tokyo 112-8551, Japan;
e-mail: katori@phys.chuo-u.ac.jp
}
and 
Hideki Tanemura
\footnote{
Department of Mathematics and Informatics,
Faculty of Science, Chiba University, 
1-33 Yayoi-cho, Inage-ku, Chiba 263-8522, Japan;
e-mail: tanemura@math.s.chiba-u.ac.jp
}}
%%%%%%%%%%%%%%%%%%%%%%%%%%%%%%%%
\date{15 January 2013}
%%%%%%%%%%%%%%%%%%%%%%%%%%%%%%%%
\pagestyle{plain}
\maketitle
\begin{abstract}
Dyson's Brownian motion model with the parameter $\beta=2$,
which we simply call the Dyson model in the present paper,
is realized as an $h$-transform of the absorbing Brownian motion
in a Weyl chamber of type A.
Depending on initial configuration
with a finite number of particles,
we define a set of entire functions
and introduce a martingale for a system of
independent complex Brownian motions (CBMs),
which is expressed by a determinant of 
a matrix with elements given by
the conformal transformations of CBMs by the entire functions.
We prove that the Dyson model can be represented
by the system of independent CBMs weighted
by this determinantal martingale.
 From this CBM representation,
the Eynard-Mehta-type correlation kernel
is derived and the Dyson model is
shown to be determinantal.
The CBM representation is a useful extension
of $h$-transform, since it works also in
infinite particle systems.
Using this representation, we 
prove the tightness of a series of processes,
which converges to the Dyson model with 
an infinite number of particles,
and the noncolliding property of the limit process.

\vskip 0.3cm
\noindent{\bf AMS 2000 Subject classifications:}
60J70, 60G44, 82C22, 32A15, 15A52 

\vskip 0.3cm
\noindent{\bf Keywords:}
the Dyson model, $h$-transform, 
complex Brownian motions, entire functions,
conformal martingales, 
determinantal process, tightness, noncolliding property
\end{abstract}

\clearpage

%%%%%%%%%%%%%%%%%%%%%%%%%%%%%%%%%%%%%%%%%%%%%%%%%%%%%%%%%%
%%%%%%%%%%%%%%%%%%%%%%%%%%%%%%%%%%%%%%%%%%%%%%%%%%%%%%%%%%
%%%  SEC1   %%%%%%%%%%%%%%%%%%%%%%%%%%%%%%%%%%%%%%%%%%%%%%
%%%%%%%%%%%%%%%%%%%%%%%%%%%%%%%%%%%%%%%%%%%%%%%%%%%%%%%%%%
\SSC{Introduction and Results}\label{chap: Introduction}%%
%%%%%%%%%%%%%%%%%%%%%%%%%%%%%%%%%%%%%%%%%%%%%%%%%%%%%%%%%%

Dyson's Brownian motion model 
is a one-parameter family of the systems 
of one-dimensional Brownian motions
with long-ranged repulsive interactions,
whose strength is represented by
a parameter $\beta >0$.
The system solves the following 
stochastic differential equations (SDEs), 
\begin{equation}
dX_i(t)=dB_i(t) 
+\frac{\beta}{2} 
\sum_{\substack{1 \leq j \leq n: \\ j \not= i}}
\frac{dt}{X_i(t)-X_j(t)}, \quad
1 \leq i \leq n, \quad
t \in [0, \infty),
\label{eqn:Dyson}
\end{equation}
where $B_i(t)$'s are independent one-dimensional
standard Brownian motions \cite{Dys62,Spo87}.
In the present paper we consider the case with $\beta=2$,
since in this special case the system is
realized by the following three processes
 \cite{KT04,KT07},
\begin{description}
\item{(i)} \quad 
the process of eigenvalues of Hermitian matrix-valued diffusion 
process in the Gaussian unitary ensemble (GUE) 
\cite{Dys62,Meh04,For10},

\item{(ii)} \quad 
the system of one-dimensional Brownian motions 
conditioned never to collide with each other \cite{Gra99}, 

\item{(iii)} \quad 
the \textit{harmonic transform} of the absorbing Brownian 
motion in a Weyl chamber of type A$_{n-1}$ \cite{Gra99}, 
\[
\mathbb{W}_n^{\textrm{A}}=\{ \mib{x} \in \mathbb{R}^n : x_1 < x_2 < \cdots < x_n \}, 
\]
with a harmonic function given by the Vandermonde determinant
\begin{equation}
h(\mib{x}) = \prod_{1 \leq i< j \leq n}(x_j-x_i)
= \det_{1 \leq i, j \leq n}\left[x_{j}^{i-1} \right].
\label{eqn:Vandermonde}
\end{equation}
\end{description}
In the family of particle systems (\ref{eqn:Dyson}), the case with $\beta=2$ 
plays the role which is similar 
to that of the three-dimensional ($\nu=1/2$) Bessel process
in the family of Bessel processes 
with parameter $\nu >-1$ \cite{KT11b}.
We call the case with $\beta=2$ of Dyson's Brownian motion model
simply \textit{the Dyson model} in this paper.

Let $\mM$ be the space of nonnegative integer-valued Radon measures 
on $\mathbb{R}$, which is a Polish space with the vague topology.
Any element $\xi$ of $\mM$ can be represented as
$\xi(\cdot) = \sum_{i \in \mathbb{I}}\delta_{x_i}(\cdot)$
with a countable index set $\mathbb{I}$ and
a sequence of points in $\mathbb{R}$, $\mib{x} =(x_i)_{i \in \mathbb{I}}$ 
satisfying $\xi(K)=\sharp\{x_i: x_i \in K \} < \infty$ 
for any compact subset $K \subset \mathbb{R}$.
In this paper the cardinality of a finite set $S$
is denoted by $\sharp S$.
We call an element $\xi$ of $\mM$ an unlabeled configuration,
and a sequence $\mib{x}$ a labeled configuration.
We write the restriction of configuration 
in $A \subset \mathbb{R}$ as
$(\xi\cap A) (\cdot)=\sum_{i \in \mathbb{I} : x_i \in  A}
\delta_{x_i}(\cdot)$, 
a shift of configuration by $u \in \mathbb{R}$ as
$\tau_u \xi(\cdot)=\sum_{i \in \mathbb{I}} \delta_{x_i+u}(\cdot)$,
and a square of configuration as 
$\xi^{\langle 2 \rangle}(\cdot)
=\sum_{i \in \mathbb{I}} \delta_{x_i^2}(\cdot)$, respectively.
The set of $\mM$-valued 
continuous functions defined on $[0,\infty)$
is denoted by $\textrm{C}([0,\infty)\to \mM)$,
which has the topology of uniform convergence on any compact sets.
For $\xi \in \mM$ with $\xi(\mathbb{R}) \in \mathbb{N} \equiv \{1,2, \dots\}$, 
we introduce a one-parameter family of \textit{entire functions}
of $z \in \mathbb{C}$ \cite{L96} parameterized by $u \in \mathbb{C}$,
$\{ \Phi_{\xi}^{u}(z) : u \in \mathbb{C}\}$, in which
\begin{equation}
\Phi_{\xi}^{u}(z)= 
\prod_{r \in\textrm{supp } \xi \cap\{u\}^\textrm{c}}
\left( 1 - \frac{z-u}{r-u} \right)^{\xi(\{r\})}
\label{eqn:entireF}
\end{equation}
with $\textrm{supp }  \xi = \{r \in \mathbb{R} : \xi(\{r\}) > 0\}$.
The zero set of the function (\ref{eqn:entireF}) is 
$\textrm{supp } \xi \cap \{u \}^{\textrm{c}}$.

With the SDEs (\ref{eqn:Dyson}), 
we consider the diffusion process 
$\Xi(t)=\sum_{i \in \mathbb{I}} \delta_{X_i(t)}$ in $\mM$
and the process under the initial configuration 
$\xi=\sum_{i \in \mathbb{I}} \delta_{x_i} \in\mM$
is denoted by $(\Xi(t), \mathbb{P}_{\xi})$.
We write the expectation with respect to $\mathbb{P}_{\xi}$ as
$\mathbb{E}_{\xi}$.
We introduce a filtration $\{ {\cal F}(t) \}_{t\in [0,\infty)}$
on the space $\textrm{C}([0,\infty)\to \mM)$
defined by ${\cal F}(t) = \sigma (\Xi(s), s\in [0,t])$.
Let $\textrm{C}_0(\mathbb{R})$ be the set of all continuous
real-valued functions with compact supports.
For any integer $M \in \mathbb{N}$,
a sequence of times
$\mib{t}=(t_1,t_2,\dots,t_M)$ with 
$0 < t_1 < \cdots < t_M \leq T < \infty$,
and a sequence of functions
$\mib{f}=(f_{t_1},f_{t_2},\dots,f_{t_M}) \in \textrm{C}_0(\mathbb{R})^M$,
the moment generating function of multitime distribution
of $(\Xi(t), \mathbb{P}_{\xi})$ is defined by
\begin{equation}
{\Psi}_{\xi}^{\mib{t}}[\mib{f}]
\equiv \mathbb{E}_{\xi} \left[\exp \left\{ \sum_{m=1}^{M} 
\int_{\mathbb{R}} f_{t_m}(x) \Xi(t_m, dx) \right\} \right].
\label{def:GF}
\end{equation}
We put 
$\mM_0= \{ \xi\in\mM : \xi(\{x\})\le 1 \mbox { for any }  x\in\mathbb{R}
\}$. 
Since any element $\xi$ of $\mM_0$ is determined uniquely 
by its support, 
it is identified with a countable subset 
$\{x_i\}_{i\in \mathbb{I}}$ of $\mathbb{R}$.
When $\xi = \sum_{i \in \mathbb{I}} \delta_{u_i} \in \mM_0$,
(\ref{eqn:entireF}) gives
\begin{equation}
\Phi_{\xi}^{u_i}(u_j)=\delta_{ij},
\quad i, j \in \mathbb{I}.
\label{eqn:Kronecker1}
\end{equation}

For a finite index set $\mathbb{I}$ and
$\mib{u}=(u_i)_{i \in \mathbb{I}}, u_i \in \mathbb{R}$,
let $Z_i(t), t\ge 0$, $i\in \mathbb{I}$ 
be a sequence of independent complex Brownian motions (CBMs)
on a probability space $(\Omega, {\cal F}, \textbf{P}_{\mib{u}})$
with $Z_i(0)=u_i$.
We write the expectation with respect to $\textbf{P}_{\mib{u}}$ as $\textbf{E}_{\mib{u}}$.
The real part and the imaginary part
of $Z_i(t)$ are denoted by
$V_i(t)={\rm Re} Z_i(t)$ and
$W_i(t)={\rm Im} Z_i(t)$,  
respectively, $i \in \mathbb{I}$,
which are independent one-dimensional
standard Brownian motions.
For any sequences $(u_i)_{i\in\mathbb{I}}$ and $x\in\mathbb{R}$,
if we set $\xi=\sum_{i \in \mathbb{I}} \delta_{u_i}$,
\begin{equation}
\Phi_{\xi}^{x}(Z_i(\cdot)), \ i \in \mathbb{I} 
\mbox{ are independent conformal local martingales,}
\label{eqn:martingale}
\end{equation}
since $\Phi_{\xi}^{x}$ is an entire function.
Each of them is a time change of a CBM \cite{RY98}.
When $\xi=\sum_{i \in \mathbb{I}} \delta_{u_i} \in \mM_0$,
combination of (\ref{eqn:Kronecker1}) and (\ref{eqn:martingale})
gives, for $0 < t \leq T < \infty$,
\begin{equation}
\textbf{E}_{\mib{u}}[\Phi_{\xi}^{u_i}(Z_j(t))]
=\textbf{E}_{\mib{u}}[\Phi_{\xi}^{u_i}(Z_j(0))]
=\Phi_{\xi}^{u_i}(u_j)=\delta_{ij},
\quad i, j \in \mathbb{I}.
\label{eqn:Kronecker2}
\end{equation}

A key observation for the present study is that the equality
\begin{equation}
\det_{1\le i,j \le \xi(\mathbb{R})}\Big[ \Phi_{\xi}^{u_i}(z_j)\Big]
= \frac{h(\mib{z})}{h(\mib{u})}
\label{eqn:determinant2}
\end{equation}
holds for any $\xi=\sum_{i=1}^{\xi(\mathbb{R})} \delta_{u_i}$
with $\xi(\mathbb{R}) \in \mathbb{N}$, $\mib{u}=(u_1, \dots, u_{\xi(\mathbb{R})}) \in \mathbb{W}_{\xi(\mathbb{R})}^{\rm A}$
and $\mib{z}=(z_1, \dots, z_{\xi(\mathbb{R})}) \in \mathbb{C}^{\xi(\mathbb{R})}$.
This is proved as follows.
Let $\xi(\mathbb{R})=n$ and 
\[
H(\mib{u},\mib{z})= \det_{1\le i,j \le n} 
\Big[ \prod_{1\le k \le n: k\not= i}(u_k-z_j) \Big],
\quad \mib{u}, \mib{z}\in \mathbb{C}^n.
\] 
Since $H$ is a polynomial function with degree $n(n-1)$ 
satisfying the conditions that $H(\mib{u},\mib{u})=(-1)^{n(n-1)/2} h(\mib{u})^2$, and 
$H(\mib{u}, \mib{z})=0, $
if $u_i=u_j$ or $z_i=z_j$ for some $i, j$ with $1\le i<j\le n$, 
we find $H(\mib{u},\mib{z})=(-1)^{n(n-1)/2} h(\mib{u}) h(\mib{z})$.
(It is a special case of the determinantal identity
given as Lemma 2.2 in \cite{Kra90}.)
Since the LHS of (\ref{eqn:determinant2})
is equal to 
$(-1)^{n(n-1)/2}H(\mib{u},\mib{z})/h(\mib{u})^2$
by definition (\ref{eqn:entireF}) for $\xi \in \mM_0$,
we obtain (\ref{eqn:determinant2}).
This equality implies that from a harmonic function $h(\cdot)$
given by (\ref{eqn:Vandermonde}),
we have a martingale for a system of 
independent CBMs $\{Z_i(\cdot) : 1 \leq i \leq {\xi(\mathbb{R})} \}$
in the determinantal form, 
$\det_{1 \leq i, j \leq {\xi(\mathbb{R})}} 
[ \Phi_{\xi}^{u_i}(Z_j(\cdot))]$.

Let $\1 (\omega)$ be
the indicator function of a condition $\omega$;
$\1(\omega)=1$ if $\omega$ is satisfied 
and $\1(\omega)=0$ otherwise,
and $\mathbb{I}_p=\{1,2,\dots,p\}$ for $p\in \mathbb{N}$. 
The main theorem of the present paper is 
the following.

%%%%%%%%%%%%%%%%%%%%%%%%%%%%%%%%%%%%%%%%%%
%%%%%%%%% Theorem %%%%%%%%%%%%%%%%%%%%%%%%%%%%
%%%%%%%%%%%%%%%%%%%%%%%%%%%%%%%%%%%%%%%%%%

\begin{theorem} \label{thm:main}
Suppose that $\xi=\sum_{i=1}^{\xi(\mathbb{R})} \delta_{u_i} 
\in \mM_0$ with $\xi(\mathbb{R}) \in \mathbb{N}$.
Let $0 < t \leq T < \infty$.
For any ${\cal F}(t)$-measurable bounded function $F$
we have
\begin{equation}
\mathbb{E}_{\xi}\left[F\big(\Xi(\cdot)\big)\right]
=\bE_{\mib{u}} \left[F \left( \sum_{i=1}^{\xi(\mathbb{R})} \delta_{V_i(\cdot)} \right)
\det_{1\le i,j \le {\xi(\mathbb{R})}}
\Big[\Phi_{\xi}^{u_i}(Z_j(T))\Big]\right].
\label{eqn:main}
\end{equation}
In particular, the moment generating function (\ref{def:GF}) 
is given by
\begin{eqnarray}
&&{\Psi}_{\xi}^{\mib{t}}[\mib{f}]
=\sum_{p=0}^{\xi(\mathbb{R})} 
\sum_{(\mathbb{J}_m )_{m=1}^M}
%\sum_{\substack{(\mathbb{J}_m )_{m=1}^M: \\
%\bigcup_{m=1}^M \mathbb{J}_m = \mathbb{I}_{p}}}
\int_{\mathbb{W}_{p}^{\rm A}} \xi^{\otimes p}(d\mib{v})
\bE_{\mib{v}} \left[
\prod_{m=1}^M \prod_{j_m\in \mathbb{J}_m} \chi_{t_m}(V_{j_m}(t_m))
\det_{i,j \in \mathbb{I}_p}
\Big[\Phi_{\xi}^{v_i}(Z_j(T))\Big]
 \right],
\nonumber\\
\label{eqn:cBMexp1}
\end{eqnarray}
where
%$\chi_{t_m}(\cdot)=e^{f_{t_m}(\cdot)}-1$, $1\le m \le M$.
$\chi_{t_m}(\cdot)=e^{f_{t_m}(\cdot)}-1$, $1\le m \le M$,
and  the second summation in the right hand side of the above equation 
runs over all $\mathbb{J}_1,\, \mathbb{J}_2, \dots,\, \mathbb{J}_M \subset \mathbb{I}_{p}$
 with $\displaystyle{\bigcup_{m=1}^M} \mathbb{J}_m = \mathbb{I}_{p}$.
\end{theorem}
%%%%%%%%%%%%%%%%%%%%%%%%%%%%%%%%%%%%%%%%%%%%%%%%%%%%%%%%%%%%

We call the above results 
\textit{the complex Brownian motion (CBM)
representations} of the Dyson model.
In order to show the simplest application of this 
representation, we consider the density function
at a single time for $(\Xi(t), \mathbb{P}_{\xi})$
denoted by $\rho_{\xi}(t, x)$.
It is defined as a continuous function of $x \in \mathbb{R}$
for $0 < t \leq T < \infty$ such that
for any $\chi \in \textrm{C}_0(\mathbb{R})$
\begin{equation}
\mathbb{E}_{\xi} \left[
\int_{\mathbb{R}} \, \chi(x) \Xi(t, dx) \right]
=\int_{\mathbb{R}} dx \, \chi(x) \rho_{\xi}(t, x).
\label{eqn:density1}
\end{equation}
By (\ref{eqn:Kronecker2}), 
the equality (\ref{eqn:main}) gives the following
expression for (\ref{eqn:density1})
\begin{eqnarray}
&& \int_{\mathbb{R}} \xi(dv)
\textbf{E}_{v} \Big[ \chi(V(t))
\Phi_{\xi}^{v}(Z(t)) \Big]
\nonumber\\
&& \qquad
= \int_{\mathbb{R}} \xi(dv) \int_{\mathbb{R}} dx \,
p_{0,t}(v, x) \int_{\mathbb{R}} dw \,
p_{0,t}(0, w) \chi(x)
\Phi_{\xi}^{v}(x+\sqrt{-1} w),
\nonumber
\end{eqnarray}
where 
$\displaystyle{p_{s,t}(x,y)=\frac{e^{-(y-x)^2/2 (t-s)}}{\sqrt{2 \pi (t-s)}}}$,
$0 \leq s < t, x, y \in \mathbb{R}$,
since $V(0)={\rm Re} Z(0)=v \in \textrm{supp } \xi$
and $W(0)={\rm Im} Z(0)=0$.
Then, if we define the function
\begin{equation}
\sG_{s,t}(x,y)= \int_{\mathbb{R}}\xi(dv)p_{0,s}(v,x)
\int_{\mathbb{R}} dw \ p_{0,t}(0,w)\Phi_{\xi}^{v}(y+\sqrt{-1}w)
\label{eqn:sGtt}
\end{equation}
for $(x,y) \in \mathbb{R}^2, (s,t) \in (0, T]^2$,
we obtain the expression for the density
function
\[
\rho_{\xi}(t, x)= \sG_{t,t}(x,x),
\quad x \in \mathbb{R}, \quad 0 < t \leq T < \infty
\]
for any initial configuration $\xi \in \mM_0$.
The above calculation will be fully generalized and
the following formula can be derived from our CBM representations.

%%%%%%%%%%%%%%%%%%%%%%%%%%%%%%%%%%%%%%%%%%%%%%%%%%
%%%%%%%%% Corollary %%%%%%%%%%%%%%%%%%%%%%%%%%%%
%%%%%%%%%%%%%%%%%%%%%%%%%%%%%%%%%%%%%%%%%%%%%%%%%%
\begin{corollary} \label{thm:2time}
Suppose that $\xi \in \mM_0$ with $\xi(\mathbb{R}) \in \mathbb{N}$.
Let
\begin{equation}
\mathbb{K}_{\xi}(s, x; t, y) = 
\sG_{s,t}(x,y)-\1(s>t)p_{t,s}(y,x).
\label{eqn:Kgen1}
\end{equation}
Then the moment generating function (\ref{eqn:cBMexp1}) for the
multitime distribution
is given by a Fredholm determinant
\begin{equation}
{\Psi}_{\xi}^{\mib{t}}[\mib{f}]
=\mathop{\Det}_
{\substack
{(s,t)\in \{t_1,t_2,\dots, t_M\}^2, \\
(x,y)\in \mathbb{R}^2}
}
 \Big[\delta_{st} \delta_x(y)
+ \mathbb{K}_{\xi}(s,x;t,y) \chi_{t}(y) \Big].
\label{eqn:Fred}
\end{equation}
\end{corollary}
%%%%%%%%%%%%%%%%%%%%%%%%%%%%%%%%%%%%%%%%%%%%%%%%%%%%%%%%%%%%

By definition of Fredholm determinant 
(see, for example, \cite{For10})
the moment generating function (\ref{eqn:Fred})
can be expanded with respect 
to $\chi_{t_m}(\cdot), 1 \leq m \leq M$, as
\begin{equation}
{\Psi}_{\xi}^{\mib{t}}[\mib{f}]
=\sum_
{\substack
{N_m \geq 0, \\ 1 \leq m \leq M} }
\int_{\prod_{m=1}^{M} \mathbb{W}^{\rm A}_{N_{m}}}
\prod_{m=1}^{M} \left\{ d \mib{x}_{N_m}^{(m)}
\prod_{i=1}^{N_{m}} 
\chi_{t_m} \Big(x_{i}^{(m)} \Big) \right\}
\rho_{\xi} 
\Big( t_{1}, \mib{x}^{(1)}_{N_1}; \dots ; t_{M}, \mib{x}^{(M)}_{N_M} \Big),
\label{eqn:ch_rho}
\end{equation}
with
\begin{equation}
\rho_{\xi} \Big(t_1,\mib{x}^{(1)}_{N_1}; \dots;t_M,\mib{x}^{(M)}_{N_M} \Big) 
=\det_{\substack
{1 \leq i \leq N_{m}, 1 \leq j \leq N_{n}, \\
1 \leq m, n \leq M}
}
\Bigg[
\mathbb{K}_{\xi}(t_m, x_{i}^{(m)}; t_n, x_{j}^{(n)} )
\Bigg],
\label{eqn:determinant}
\end{equation}
where $\mib{x}^{(m)}_{N_m}$ denotes
$(x^{(m)}_1, \dots, x^{(m)}_{N_m})$
and $d\mib{x}^{(m)}_{N_m}= \prod_{i=1}^{N_m} dx^{(m)}_i$,
$1 \leq m \leq M$.
The functions $\rho_{\xi}$'s are multitime 
correlation functions, and ${\Psi}_{\xi}^{\mib{t}}[\mib{f}]$ 
can be regarded as 
a generating function of them.
The function $\mathbb{K}_{\xi}$ 
given by (\ref{eqn:Kgen1}) with (\ref{eqn:sGtt}) 
is thus called the \textit{correlation kernel} \cite{KT10}.
In general, when the moment generating function for
the multitime distribution is given by
a Fredholm determinant,
the process is said to be \textit{determinantal} \cite{KT07,KT10}.
The results by Eynard and Mehta reported in \cite{EM98}
for a multi-layer matrix model can be regarded as the theorem that
the Dyson model is determinantal for the special
initial configuration $\xi =\xi(\mathbb{R}) \delta_0$, \textit{ i.e.},
all particles are put at the origin, 
for any $\xi(\mathbb{R}) \in \mathbb{N}$.
The correlation kernel is expressed by using
the Hermite orthogonal polynomials \cite{NF98}.
The present authors proved that, 
for any fixed initial configuration $\xi \in \mM$
with $\xi(\mathbb{R}) \in \mathbb{N}$,
the Dyson model $(\Xi(t), \mathbb{P}_{\xi})$
is determinantal,
in which the correlation kernel is given by
\begin{eqnarray}
\mathbb{K}_{\xi}(s,x; t,y)
&=& \frac{1}{2 \pi \sqrt{-1}} 
\oint_{\Gamma(\xi)} dz \, p_{0,s}(z, x)
\int_{\mathbb{R}} d w \, p_{0,t}(w, -\sqrt{-1} y)
\nonumber\\
&\times&
\frac{1}{\sqrt{-1} w-z} 
\prod_{r \in \textrm{supp } \xi} 
\left( 1- \frac{\sqrt{-1} w-z}{r-z} \right)
-\1(s>t) p_{t,s}(y,x),
\label{eqn:KgenC}
\end{eqnarray}
where $\Gamma(\xi)$ is a closed contour on the
complex plane $\mathbb{C}$ encircling the points in 
$\textrm{supp } \xi$ on the real line $\mathbb{R}$
once in the positive direction
(Proposition 2.1 in \cite{KT10}).
In order to derive (\ref{eqn:KgenC}), we used the
integral representations of multiple
Hermite polynomials given by 
Bleher and Kuijlaars \cite{BK05}.

In the present paper, we assume $\xi \in \mM_0$
preventing the initial configuration
from having any multiple points.
This restriction is only for simplicity of calculation.
(Note that, if $\xi \in \mM_0$,
the Cauchy integrals in (\ref{eqn:KgenC}) can be
readily performed and the expression 
(\ref{eqn:Kgen1}) with (\ref{eqn:sGtt}) is obtained.)
The fact that we would like to report here
is that, although the theory of (multiple-)orthogonal functions
are very useful to analyze determinantal processes
\cite{KT10,KT09,KT11a},
it is not necessary to deriving
the Eynard-Mehta-type determinantal expressions
for multitime correlation functions.
The essential point may be the extension of
$h$-transform to the conformal 
martingale of CBMs in the determinantal form
$\det_{1 \leq i, j \leq \xi(\mathbb{R})} [\Phi_{\xi}^{u_i}(Z_j(\cdot))]$,
which we have named the CBM representation.
In other words, the proof of Corollary \ref{thm:2time}
will provide a probability-theoretical derivation
of the \textit{Eynard-Mehta-type correlation kernel}.

We gave useful sufficient conditions of $\xi$ 
in \cite{KT10} so that the Dyson model is well defined
as a determinantal process even if $\xi(\mathbb{R})=\infty$.
For $L>0, \alpha>0$ and $\xi\in\mM$ we put
\[
M(\xi, L)=\int_{[-L,L]\setminus\{0\}} \frac{\xi(dx)}{x},
\qquad
M_\alpha(\xi, L)
=\left( \int_{[-L,L]\setminus\{0\}} 
\frac{\xi(dx)}{|x|^\alpha}\right)^{1/\alpha},
\]
and
$\displaystyle{M(\xi) = \lim_{L\to\infty}M(\xi, L)}$,
$\displaystyle{M_\alpha(\xi)= \lim_{L\to\infty}M_\alpha(\xi, L)}$,
if the limits finitely exist. Then
%%%%%%%%%%%%%%%%%%%%%%%%%%%%%%%%%%%%%%%%%%%%%%%%%%%%%%%%%%%%%
%%%%%%%%%%%%%%Conditions %%%%%%%%%%%%%%%%%%%%%%%%%%%%%%%%%%%%
%%%%%%%%%%%%%%%%%%%%%%%%%%%%%%%%%%%%%%%%%%%%%%%%%%%%%%%%%%%%%
\vskip 3mm
%%%%%%%%%%%%%%%%%%%%%%%%%%%%%%%%%%%%%%%%%%%%%%%%%%%%%%%%%%%%%
\noindent (\textbf{C.1})
there exists $C_0 > 0$ such that
$|M(\xi,L)|  < C_0$, $L>0$,

\vskip 3mm

\noindent (\textbf{C.2}) (i) 
there exist $\alpha\in (1,2)$ and $C_1>0$ such that
$
M_\alpha(\xi) \le C_1,
$ \\
\noindent (ii) 
there exist $\beta >0$ and $C_2 >0$ such that
\[
M_1(\tau_{-a^2} \xi^{\langle 2 \rangle}) \leq C_2
(\max\{|a|, 1\})^{-\beta}
\quad \forall a \in \textrm{supp } \xi.
\]
%%%%%%%%%%%%%%%%%%%%%%%%%%%%%%%%%%%%%%%%%%%%%%%%%%%%%%%%%%%%%
\vskip 3mm
%%%%%%%%%%%%%%%%%%%%%%%%%%%%%%%%%%%%%%%%%%%%%%%%%%%%%%%%%%%%%
%%%%%%%%%%%%%%%%%%%%%%%Conditions%%%%%%%%%%%%%%%%%%%%%%%%%%%%
%%%%%%%%%%%%%%%%%%%%%%%%%%%%%%%%%%%%%%%%%%%%%%%%%%%%%%%%%%%%

\noindent
It was shown that, 
if $\xi\in\mM_0$ satisfies the 
conditions $(\textbf{C.1})$ and $(\textbf{C.2})$, 
then for $a\in\mathbb{R}$ and $z\in\mathbb{C}$, 
$\displaystyle{\Phi_{\xi}^{a}(z)\equiv \lim_{L\to\infty}
\Phi_{\xi \cap [a-L, a+L]}^{a}(z)}$
finitely exists,
and
\begin{equation}
|\Phi_{\xi}^{a}(z)|\le C \exp \bigg\{c(|a|^\theta +|z|^\theta) \bigg\}
\left| \frac{z}{a}\right|^{\xi(\{0\})} \left|\frac{a}{a-z} \right|,
\quad a\in \textrm{supp } \xi, \ z \in\mathbb{C},
\label{estimate:Phi}
\end{equation}
for some $c, C>0$ and $\theta\in 
(\max\{\alpha, (2-\beta)\},2)$, 
which are determined by the constants $C_0, C_1, C_2$
and the indices $\alpha, \beta$ 
in the conditions (Lemma 4.4 in \cite{KT10}).
We have noted that in the case that $\xi \in \mM_0$ 
satisfies the conditions 
$(\textbf{C.1})$ and $(\textbf{C.2})$
with constants $C_0, C_1, C_2$ and indices $\alpha$ and $\beta$,
then $\xi \cap [-L, L], \forall L > 0$ does as well.
Hence we can obtain the convergence of moment generating functions
${\Psi}_{\xi \cap [-L, L]}^{\mib{t}}[\mib{f}] \to
{\Psi}_{\xi}^{\mib{t}}[\mib{f}]$ as $L \to \infty$, 
which implies the convergence of
the probability measures $\mathbb{P}_{\xi\cap [-L, L]}
\to \mathbb{P}_{\xi}$ in $L \to \infty$
in the sense of finite dimensional distributions.
Moreover,
even if $\xi(\mathbb{R})=\infty$, $\mathbb{K}_{\xi}$ given by 
(\ref{eqn:Kgen1}) with (\ref{eqn:sGtt}) 
is well-defined as a correlation kernel and
dynamics of the Dyson model with an infinite number of
particles $(\Xi(t), \mathbb{P}_{\xi})$ exists 
as a determinantal process \cite{KT10}.

The CBM representation is indeed a non-trivial
extension of $h$-transforms,
since it works also in infinite particle systems.

%%%%%%%%%%%%%%%%%%%%%%%%%%%%%%%%%%%%%%%%%%%%%%%%%%%%%%%%%%%%%
%%%%%%Cor %%%%%%%%% Cor infinite cBM rep %%%%%%%%%%%%%%
%%%%%%%%%%%%%%%%%%%%%%%%%%%%%%%%%%%%%%%%%%%%%%%%%%%%%%%%%%%%%
\begin{corollary}\label{cor:infinite_cBM}
Suppose that the initial configuration
$\xi\in \mM_0$ satisfies the conditions 
$({\bf C.1})$ and $({\bf C.2})$.
Then the expression (\ref{eqn:main}) is valid also
in the case with $\xi(\mathbb{R})=\infty$, 
if $F$ is represented as
\[
F(\Xi(\cdot))= G \left(\int_{\mathbb{R}}\phi_1(x)\Xi(t_1, dx),
\int_{\mathbb{R}}\phi_2(x)\Xi(t_2, dx),\dots, 
\int_{\mathbb{R}}\phi_k(x)\Xi(t_k, dx) \right)
\]
with $\phi_i \in C_0(\mathbb{R})$, $1\le i \le k$ 
and a polynomial function $G$ on $\mathbb{R}^k$, $k\in\mathbb{N}$.
\end{corollary}
%%%%%%%%%%%%%%%%%%%%%%%%%%%%%%%%%%%%%%%%%%%%%%%%%%%%%%%%%%%%

In order to demonstrate the usefulness of the CBM
representations 
to characterize infinite particle systems, 
we show that the following estimate
is readily obtained from the expression (\ref{eqn:main}).
Let $\textrm{C}_{0}^{\infty}(\mathbb{R})$ be the set of all infinitely 
differentiable real functions with compact supports.

%%%%%%%%%%%%%%%%%%%%%%%%%%%%%%%%%%%%%%%%%%%%%%%%%%%%%%%%%%%%%
%%%%%%PROP %%%%%%%%% Prop 4-th moment estimate %%%%%%%%%%%%%%
%%%%%%%%%%%%%%%%%%%%%%%%%%%%%%%%%%%%%%%%%%%%%%%%%%%%%%%%%%%%%
\begin{proposition}\label{prop:4-th_moment}
Suppose that the initial configuration
$\xi\in \mM_0$ satisfies the conditions 
$({\bf C.1})$ and $({\bf C.2})$
with constants $C_0, C_1, C_2$ and indices $\alpha, \beta$.
Then for any $T>0$ and $\varphi \in \textrm{C}_0^{\infty}(\mathbb{R})$
there exists a positive constant $C=C(C_0,C_1,C_2,\alpha,\beta,T,\varphi)$,
which is independent of $s, t$, such that
\begin{equation}
\mathbb{E}_{\xi} \Bigg[
\left|
\int_{\mathbb{R}} \varphi(x) \Xi(t, dx)
-\int_{\mathbb{R}} \varphi(x) \Xi(s, dx)
\right|^4 
\Bigg]
\leq C |t-s|^2,
\quad \forall s, t \in [0,T].
\label{eqn:KolCri}
\end{equation}
\end{proposition}
%%%%%%%%%%%%%%%%%%%%%%%%%%%%%%%%%%%%%%%%%%%%%%%%%%%%%%%%%%%%%
%%%%%%PROP %%%%%%%%% Prop 4-th moment estimate %%%%%%%%%%%%%%
%%%%%%%%%%%%%%%%%%%%%%%%%%%%%%%%%%%%%%%%%%%%%%%%%%%%%%%%%%%%%

By a criterion of Kolmogorov 
(see, for example, \cite{Kal97,Bil99}), 
Proposition \ref{prop:4-th_moment} implies 
that the sequence of the process 
$(\Xi(t),\mathbb{P}_{\xi\cap [-L, L]})$, $L \in \mathbb{N}$ is tight
in $\textrm{C}([0,\infty)\to \mM)$.
Then we can conclude the following.

%%%%%%%%%%%%%%%%%%%%%%%%%%%%%%%%%%%%%%%%%%%%%%%%%%%%%%%%%%%%
%%%%%%%%% THEOREM continuity %%%%%%%%%%%%%%%%%%%%%%%%%%%%%%%
%%%%%%%%%%%%%%%%%%%%%%%%%%%%%%%%%%%%%%%%%%%%%%%%%%%%%%%%%%%%
\begin{theorem}\label{Theorem:Tightness}
Suppose that $\xi\in \mM_0$ satisfies the conditions 
$({\bf C.1})$ and $({\bf C.2})$.
Then the process $(\Xi(t),\mathbb{P}_{\xi\cap [-L, L]})$
converges to the process $(\Xi(t),\mathbb{P}_\xi)$
weakly on $\textrm{C}([0, \infty)\to \mM)$ as $L \to\infty$.
In particular, the process $(\Xi(t),\mathbb{P}_\xi)$
has a modification which is almost-surely continuous on $[0,\infty)$
with $\Xi(0)=\xi$.
\end{theorem}
%%%%%%%%%%%%%%%%%%%%%%%%%%%%%%%%%%%%%%%%%%%%%%%%%%%%%%%%%%%%
%%%%%%%%% THEOREM continuity %%%%%%%%%%%%%%%%%%%%%%%%%%%%%%%
%%%%%%%%%%%%%%%%%%%%%%%%%%%%%%%%%%%%%%%%%%%%%%%%%%%%%%%%%%%%

Finally in the present paper, we show that
the noncolliding property of the Dyson model 
with an infinite number of particles is obtained 
by using the CBM representations.

%%%%%%%%%%%%%%%%%%%%%%%%%%%%%%%%%%%%%%%%%%%%%%%%%%%%%%%%%%%%%
%%%%%%PROP %%%%%%%%% Noncolliding property %%%%%%%%%%%%%%
%%%%%%%%%%%%%%%%%%%%%%%%%%%%%%%%%%%%%%%%%%%%%%%%%%%%%%%%%%%%%
\begin{proposition}\label{prop:noncolliding}
Suppose that the initial configuration
$\xi\in \mM_0$ satisfies the conditions 
$({\bf C.1})$ and $({\bf C.2})$.
Then $\mathbb{P}_{\xi} \big[\Xi(t) \in \mM_0, \, t>0 \big]=1.$
\end{proposition}
%%%%%%%%%%%%%%%%%%%%%%%%%%%%%%%%%%%%%%%%%%%%%%%%%%%%%%%%%%%%%
%%%%%%PROP %%%%%%%%% Noncolliding property %%%%%%%%%%%%%%
%%%%%%%%%%%%%%%%%%%%%%%%%%%%%%%%%%%%%%%%%%%%%%%%%%%%%%%%%%%%%

In the following five sections,
we give the proofs of Theorem \ref{thm:main}, 
Corollary \ref{thm:2time},
Corollary \ref{cor:infinite_cBM},  
Proposition \ref{prop:4-th_moment}, and
Proposition \ref{prop:noncolliding},
respectively. 

%%%%%%%%%%%%%%%%%%%%%%%%%%%%%%%%%%%%%%%%%%%%%%%%%%%%%%%%%%
%%%%%%%%%%%%%%%%%%%%%%%%%%%%%%%%%%%%%%%%%%%%%%%%%%%%%%%%%%
%%%%%%%%%%%%%%%%%%%%%%%%%%%%%%%%%%%%%%%%%%%%%%%%%%%%%%%%%%
%%%  SEC 2   %%%%%%%%%%%%%%%%%%%%%%%%%%%%%%%%%%%%%%%%%%%%%
%%%%%%%%%%%%%%%%%%%%%%%%%%%%%%%%%%%%%%%%%%%%%%%%%%%%%%%%%%
\SSC{Proof of Theorem \ref{thm:main}}%%%%%%%%%%%%%%%%%%%%%
%%%%%%%%%%%%%%%%%%%%%%%%%%%%%%%%%%%%%%%%%%%%%%%%%%%%%%%%%%

For the proof of Theorem \ref{thm:main},
it is sufficient 
to consider the case that $F$ is given as
$F\left( \Xi(\cdot)
%\sum_{i=1}^{\xi(\mathbb{R})} \delta_{X_i(\cdot)} 
\right)
= \prod_{i=1}^M g_i(\mib{X}(t_i))$
for $M \in \mathbb{N}$, $0<t_1< \cdots <t_M \leq T < \infty$, 
with symmetric bounded measurable functions $g_i$ on $\mathbb{R}^{\xi(\mathbb{R})}$,
$1 \leq i \leq M$.
We give the proof for the case with $M=2$,
\textit{ i.e.}, for 
$\xi=\sum_{i=1}^{\xi(\mathbb{R})} \delta_{u_i},
\mib{u}=(u_1, \dots, u_{\xi(\mathbb{R})}) \in \mathbb{W}_{\xi(\mathbb{R})}^{\rm A}$
\begin{equation}
\mathbb{E}_{\xi}\left[ g_1(\mib{X}(t_1))g_2(\mib{X}(t_2))\right]
=\textbf{E}_{\mib{u}} \left[ 
g_1(\mib{V}(t_1))g_2(\mib{V}(t_2))
\det_{1\le i,j \le \xi(\mathbb{R})} \left[ \Phi_\xi^{u_i} (Z_j(T)) \right]
\right].
\label{eqn:2time-cor}
\end{equation}
The generalization for $M > 2$ is straightforward.

We use the fact that the Dyson model is obtained as an $h$-transform
of the absorbing Brownian motion in the Weyl chamber
$\mathbb{W}_{\xi(\mathbb{R})}^{\rm A}$ \cite{Gra99}.
Put $\tau= \inf\{ t>0 : \mib{V}(t)\notin \mathbb{W}_{\xi(\mathbb{R})}^{\rm A} \}$,
then the LHS of (\ref{eqn:2time-cor}) is given by
\begin{equation}
\textbf{E}_{\mib{u}} \left[
\1(\tau >t_2)g_1(\mib{V}(t_1))g_2(\mib{V}(t_2))
\frac{h(\mib{V}(t_2))}{h(\mib{u})}
\right].
\label{eqn:absBM}
\end{equation}
For a finite set $S$, we write the collection of all permutations
of elements in $S$ as $\cS(S)$.
In particular, we express $\cS(\mathbb{I}_{p})$ simply by $\cS_{p},
p \in \mathbb{N}$.
We put $\sigma(\mib{u})=(u_{\sigma(1)}, \dots, u_{\sigma(\xi(\mathbb{R}))})$
for each permutation $\sigma \in \cS_{\xi(\mathbb{R})}$.
We introduce the stopping times
\begin{equation}
\tau_{ij}=\inf\{t> 0 : V_i(t)= V_j(t)\}, \quad
1\le i < j \le \xi(\mathbb{R}).
\label{eqn:tauij}
\end{equation}
Let $\sigma_{ij} \in \cS_{\xi(\mathbb{R})}$ be the permutation of $(i,j)$.
Note that if $\sigma_{ij}(\mib{u})= \mib{u}$, 
the processes $\mib{V}(t)$ and $\sigma_{ij}(\mib{V}(t))$ are identical 
in distribution under the probability measure $\textbf{P}_{\mib{u}}$.
Then by the strong Markov property of the process $\mib{V}(t)$
and by the fact that $h$ is anti-symmetric and $g_1,g_2$ are symmetric, 
\[
\textbf{E}_{\mib{u}} \Big[
\1(\tau=\tau_{ij} \le t_2)g_1(\mib{V}(t_1))g_2(\mib{V}(t_2))
\frac{h(\mib{V}(t_2))}{h(\mib{u})}
\Big]=0.
\]
Since $\textbf{P}_{\mib{u}}(\tau_{ij}= \tau_{i' j'})=0$ if $(i,j)\not=(i',j')$,
and $\displaystyle{\tau= \min_{1\leq i<j\le \xi(\mathbb{R})}\tau_{ij}}$, 
\[
\textbf{E}_{\mib{u}} \Big[
\1(\tau \le t_2)g_1(\mib{V}(t_1))g_2(\mib{V}(t_2))
\frac{h(\mib{V}(t_2))}{h(\mib{u})} \Big]=0.
\]
Hence, (\ref{eqn:absBM}) equals
\begin{equation}
\textbf{E}_{\mib{u}} \left[g_1(\mib{V}(t_1))g_2(\mib{V}(t_2))
\frac{h(\mib{V}(t_2))}{h(\mib{u})}
\right].
\label{eqn:KM1}
\end{equation}
Then we use the equality (\ref{eqn:determinant2}) 
in (\ref{eqn:KM1}). 
Note that $V_i (t), 1 \le i \le \xi(\mathbb{R})$ and $W_i (t), 
1 \le i \le \xi(\mathbb{R})$ 
are independent. 
We can regard the probability space 
$(\Omega, {\cal F}, \textbf{P}_{\mib{v}})$ as a product of two probability spaces
$(\Omega_1, {\cal F}_1, {\rm P}_1)$ and $(\Omega_2, {\cal F}_2, {\rm P}_2)$,
and $V_i (t), 1 \le i \le \xi(\mathbb{R})$ are ${\cal F}_1$-measurable and 
$W_i (t), 1 \le i \le \xi(\mathbb{R})$ are ${\cal F}_2$-measurable.
We write ${\rm E}_{\alpha}$ for the expectation 
with respect to ${\rm P}_{\alpha}$, $\alpha=1,2$.
We see
\begin{equation}
{\rm E}_2\big[ h(\mib{Z}(t))\big] 
={\rm E}_2 \left[\det_{1\le i,j\le \xi(\mathbb{R})}
\left[ Z_j(t)^{i-1} \right] \right]
=\det_{1\le i,j\le \xi(\mathbb{R})} \left[
{\rm E}_2 \left[Z_j(t)^{i-1} \right]
\right],
\nonumber
\end{equation}
where the independence of $Z_j(t)$, $1\le j\le \xi(\mathbb{R})$,
is used in the last equality.
By binomial expansion, 
${\rm E}_2 [Z_j(t)^{i-1} ]=G(V_j (t))$ with
$
G(x) = \sum_{p=0}^{i-1}\binom{i-1}{ p}
{\rm E}_2 \Big[(\sqrt{-1}W_j (t))^{i-1-p} \Big] x^{p}.
$
Since $G(x)$ is a monic polynomial with degree $i-1$,
$
{\rm E}_2\big[ h(\mib{Z}(t))\big] = h(\mib{V}(t)).
$
Combining the above results
and the fact (\ref{eqn:martingale}), 
we have (\ref{eqn:2time-cor}).

For the proof of (\ref{eqn:cBMexp1}) with $M=2$,
we first prove that for any $N_1,N_2 \in \mathbb{N}$
\begin{eqnarray}
&&
\sum_{\substack{\mathbb{J}_1, \, \mathbb{J}_2 \subset \mathbb{I}_{\xi(\mathbb{R})}: 
\\ \sharp \mathbb{J}_1=N_1, \sharp \mathbb{J}_2=N_2}}
\mathbb{E}_{\xi}\left[
\prod_{m=1}^2 \prod_{j_m\in \mathbb{J}_m} \chi_{t_m}(X_{j_m}(t_m)) 
\right]
\nonumber\\
&&
=
\sum_{p=1}^{N_1+N_2}
\sum_{\substack{\mathbb{J}_1, \, \mathbb{J}_2 \subset \mathbb{I}_{p}: \\ 
\sharp \mathbb{J}_1=N_1, \sharp \mathbb{J}_2=N_2 \\
\mathbb{J}_1\cup \mathbb{J}_2= \mathbb{I}_{p}}}
\int_{\mathbb{W}_{p}^{\rm A}} \ 
\xi^{\otimes p}(d\mib{v}) 
\textbf{E}_{\mib{v}} \left[
\prod_{m=1}^2 \prod_{j_m\in \mathbb{J}_m} \chi_{t_m}(V_{j_m}(t_m))
\det_{i,j \in \mathbb{I}_p}
\Big[\Phi_{\xi}^{v_i}(Z_j(T))\Big] \right].
\nonumber\\
\label{eqn:cor_2time}
\end{eqnarray}
Applying (\ref{eqn:2time-cor}) with 
$\displaystyle{
g_m(\mib{x})=\sum_{\mathbb{J}_m\subset \mathbb{I}_{\xi(R)}: 
\sharp \mathbb{J}_m=N_m}\prod_{j_m \in \mathbb{J}_m}\chi_{t_m}(x_{j_m})
}$, $m=1,2$,
we see that the LHS of (\ref{eqn:cor_2time})
equals 
\begin{eqnarray}
&&\sum_{\substack{\mathbb{J}_1, \mathbb{J}_2 \subset \mathbb{I}_{\xi(\mathbb{R})}: 
\\ \sharp \mathbb{J}_1=N_1, \sharp \mathbb{J}_2=N_2}}
\textbf{E}_{\mib{u}} \left[
\prod_{m=1}^2\prod_{j_m\in \mathbb{J}_m} \chi_{t_m}(V_{j_m}(t_m))
\det_{i,j \in \mathbb{I}_{\xi(\mathbb{R})}}
\Big[\Phi_{\xi}^{u_i}(Z_j(T))\Big] \right]
\nonumber\\
&&=\sum_{p=1}^{N_1+N_2}
\sum_{\substack{\mathbb{J}_1, \, \mathbb{J}_2 \subset \mathbb{I}_{\xi(\mathbb{R})}: \\ 
\sharp (\mathbb{J}_1 \cup \mathbb{J}_2)=p, \\
\sharp \mathbb{J}_1=N_1, \sharp \mathbb{J}_2=N_2}}
\textbf{E}_{\mib{u}} \left[
\prod_{m=1}^2\prod_{j_m\in \mathbb{J}_m} \chi_{t_m}(V_{j_m}(t_m))
\det_{i,j \in \mathbb{J}_1\cup \mathbb{J}_2}
\Big[\Phi_{\xi}^{u_i}(Z_j(T))\Big] \right],
\nonumber
\end{eqnarray}
where we have used (\ref{eqn:Kronecker2}).
We see the RHS of the last equation coincides with
the RHS of (\ref{eqn:cor_2time}).
By using relation
\[
\exp\left\{ \sum_{m=1}^{2} 
\sum_{j_m=1}^{\xi(\mathbb{R})} f_{t_m}(x_{j_m}) \right\}
= \prod_{m=1}^2 \prod_{j_m=1}^{\xi(\mathbb{R})}
\Big\{\chi_{t_m}(x_{j_m})+1 \Big\}
=\sum_{\mathbb{J}_1, \, \mathbb{J}_2 \subset \mathbb{I}_{\xi(\mathbb{R})}} 
\prod_{m=1}^2\prod_{j_m \in \mathbb{J}_m} 
\chi_{t_m}(x_{j_m}),
\]
the equality (\ref{eqn:cBMexp1}) with $M=2$ is 
readily derived form (\ref{eqn:cor_2time}).
By the similar argument,
(\ref{eqn:cBMexp1}) is concluded from
(\ref{eqn:main}) for any $M >2$.
\hfill $\Box$

%%%%%%%%%%%%%%%%%%%%%%%%%%%%%%%%%%%%%%%%%%%%%%%%%%%%%%%%%%
%%%%%%%%%%%%%%%%%%%%%%%%%%%%%%%%%%%%%%%%%%%%%%%%%%%%%%%%%%
%%%  SEC3   %%%%%%%%%%%%%%%%%%%%%%%%%%%%%%%%%%%%%%%%%%%%%%
%%%%%%%%%%%%%%%%%%%%%%%%%%%%%%%%%%%%%%%%%%%%%%%%%%%%%%%%%%
\SSC{Proof of Corollary \ref{thm:2time}}%%%%%%%%%%%%%%%%%%
%%%%%%%%%%%%%%%%%%%%%%%%%%%%%%%%%%%%%%%%%%%%%%%%%%%%%%%%%%

Since the Fredholm determinant (\ref{eqn:Fred}) is explicitly 
given by (\ref{eqn:ch_rho}) with (\ref{eqn:determinant}),
(\ref{eqn:cBMexp1}) in Theorem \ref{thm:main} implies 
that, for proof of Corollary \ref{thm:2time},
it is enough to show that the following equality is established
for any $M \in \mathbb{N}, (N_1, \dots, N_M) \in \mathbb{N}^M$
\begin{eqnarray}
&& \int_{\prod_{m=1}^{M} \mathbb{W}^{\rm A}_{N_m}} 
\prod_{m=1}^{M} \left\{ d \mib{x}_{N_m}^{(m)}
\prod_{i=1}^{N_m} \chi_{t_m} 
\Big(x_{i}^{(m)} \Big) \right\}
\det_{\substack{1 \leq i \leq N_{m}, 1 \leq j \leq N_{n},
\\ 1\le m, n \le M}}
\Bigg[
\mathbb{K}_{\xi}(t_m, x_{i}^{(m)}; t_{n}, x_{j}^{(n)} )
\Bigg]
\nonumber\\
&=& \sum_{p=1}^{N}
\sum_{\substack{\sharp \mathbb{J}_m=N_m,\\
1 \leq m \leq M: \\
\bigcup_{m=1}^{M} \mathbb{J}_m= \mathbb{I}_{p}}}
\int_{\mathbb{W}_{p}^{\rm A}} \ 
\xi^{\otimes p}(d\mib{v})
\textbf{E}_{\mib{v}} \left[
\prod_{m=1}^{M} 
\prod_{j_m\in \mathbb{J}_m} \chi_{t_m}(V_{j_m}(t_m))
\det_{i,j \in \mathbb{I}_p}
\Big[\Phi_{\xi}^{v_i}(Z_j(T))\Big] \right].
\nonumber\\
\label{eqn:E=det}
\end{eqnarray}
If we take the summation of (\ref{eqn:E=det}) over
all $0 \leq N_m \leq \xi(\mathbb{R}), 1 \leq m \leq M$,
the LHS gives (\ref{eqn:ch_rho}) with (\ref{eqn:determinant})
and the RHS does (\ref{eqn:cBMexp1}).
In this section we will prove (\ref{eqn:E=det}). 
So in the following, 
we fix $M \in \mathbb{N}$, $(N_1, \dots, N_M) \in \mathbb{N}^M$.

Let $\mathbb{I}^{(1)}=\mathbb{I}_{N_1}$ and 
$\mathbb{I}^{(m)}=\mathbb{I}_{\sum_{k=1}^{m} N_{k}} 
\setminus \mathbb{I}_{\sum_{k=1}^{m-1} N_{k}}, 2 \leq m \leq M$.
Put $N=\sum_{m=1}^{M} N_m$ and
$\tau_i=\sum_{m=1}^{M} t_m \1(i \in \mathbb{I}^{(m)}),
1 \leq i \leq N$.
Then the integrand in the LHS of (\ref{eqn:E=det})
is simply written as
$\displaystyle{\prod_{i=1}^N \chi_{\tau_i}(x_i)}$
$\displaystyle{\det_{1 \leq i, j \leq N} [ \mathbb{K}_{\xi}(\tau_i, x_i; \tau_j, x_j)]}$,
and the integral 
$\int_{\prod_{m=1}^{M} \mathbb{W}^{\rm A}_{N_m}} 
\prod_{m=1}^M d \mib{x}_{N_m}^{(m)} (\cdot)$
can be replaced by
$\{\prod_{m=1}^{M} N_m ! \}^{-1}$
$\int_{\mathbb{R}^N} d \mib{x} \, (\cdot)$.
The determinant is defined using the notion of
permutations and we note that any permutation
$\sigma \in \cS_N$ can be decomposed
into a product of 
cycles.
Let the number of cycles in the decomposition
be $\ell(\sigma)$ and express $\sigma$ by
$\sigma = \sc_1 \sc_2 \cdots \sc_{\ell(\sigma)}$,
where $\sc_{\lambda}$
denotes a cyclic permutation
$\sc_{\lambda}= 
(c_\lambda(1) c_\lambda(2) \cdots c_\lambda(q_{\lambda}) ),
1 \leq q_{\lambda} \leq N, 1 \leq \lambda \leq \ell(\sigma)$.
For each $1 \leq \lambda \leq \ell(\sigma)$,
we write the set of entries 
$\{c_{\lambda}(i)\}_{i=1}^{q_{\lambda}}$ of $\sc_{\lambda}$
simply as $\{\sc_{\lambda}\}$,
in which the periodicity
$c_{\lambda}(i+q_{\lambda})=c_{\lambda}(i), 1 \leq i \leq q_{\lambda}$ 
is assumed.
By definition, for each $1 \leq \lambda \leq \ell(\sigma)$,
$c_{\lambda}(i), 1 \leq i \leq q_{\lambda}$
are distinct indices chosen from $\mathbb{I}_N$, 
$\{\sc_{\lambda}\} \cap \{\sc_{\lambda'}\} = \emptyset$
for $1 \leq \lambda \not= \lambda' \leq \ell(\sigma)$, and
$\sum_{\lambda=1}^{\ell(\sigma)} q_{\lambda}=N$.
The determinant $\displaystyle{\det_{1 \leq i,j \leq N}
[\mathbb{K}_{\xi}(\tau_i, x_i; \tau_j, x_j)]}$
is written as
\begin{eqnarray}
&& \sum_{\sigma \in \cS_N} (-1)^{N-\ell(\sigma)}
\prod_{\lambda=1}^{\ell(\sigma)}
\prod_{i=1}^{q_{\lambda}}
\mathbb{K}_{\xi}(\tau_{c_{\lambda}(i)}, x_{c_{\lambda}(i)};
\tau_{c_{\lambda}(i+1)}, x_{c_{\lambda}(i+1)})
\nonumber\\
&=& 
\sum_{\sigma \in \cS_N} (-1)^{N-\ell(\sigma)}
\prod_{\lambda=1}^{\ell(\sigma)}
\prod_{i=1}^{q_{\lambda}}
\Big\{
\sG_{\tau_{c_{\lambda}(i)},\tau_{c_{\lambda}(i+1)}}
(x_{c_{\lambda}(i)}, x_{c_{\lambda}(i+1)})
\nonumber\\
&& \qquad \qquad 
-\1(\tau_{c_{\lambda}(i)} > \tau_{c_{\lambda}(i+1)})
p_{\tau_{c_{\lambda}(i+1)},\tau_{c_{\lambda}(i)}}
(x_{c_{\lambda}(i+1)}, x_{c_{\lambda}(i)}) 
\Big\},
\label{eqn:star2}
\end{eqnarray}
where the definition (\ref{eqn:Kgen1}) of the correlation kernel
$\mathbb{K}_{\xi}$ is used.
In order to express binomial expansions for (\ref{eqn:star2}),
we introduce the following notations: 
For each cyclic permutation $\sc_{\lambda}$,
we consider a subset of $\{\sc_{\lambda}\}$,
$
\cC(\sc_{\lambda}) =
\Big\{ c_{\lambda}(i) \in \{\sc_{\lambda}\} : 
\tau_{c_{\lambda}(i)} > \tau_{c_{\lambda}(i+1)} \Big\}.
$
Choose $\textbf{M}_{\lambda}$ such that 
$\{\sc_{\lambda}\} \setminus \cC(\sc_{\lambda})
\subset \textbf{M}_{\lambda} \subset \{\sc_{\lambda}\}$,
and define $\textbf{M}_{\lambda}^{\rm c}=\{\sc_{\lambda}\} 
\setminus \textbf{M}_{\lambda}$.
Therefore if we put
\begin{eqnarray}
G(\sc_{\lambda}, \textbf{M}_{\lambda}) &=& 
\int_{\mathbb{R}^{\{\sc_{\lambda}\}}} \prod_{i=1}^{q_{\lambda}} \ \bigg\{
dx_{c_{\lambda}(i)}  \  \chi_{{\tau}_{c_{\lambda}(i)}}(x_{c_{\lambda}(i)})
p_{{\tau}_{c_{\lambda}(i+1)}, {\tau}_{c_{\lambda}(i)}}
(x_{c_{\lambda}(i+1)},x_{c_{\lambda}(i)})
^{\1(c_{\lambda}(i) \in \textbf{M}_{\lambda}^{\rm c})}
\nonumber\\
&& \qquad \qquad  \times
\sG_{{\tau}_{c_{\lambda}(i)},{\tau}_{c_{\lambda}(i+1)}}
(x_{c_{\lambda}(i)},x_{c_{\lambda}(i+1)})
^{\1(c_{\lambda}(i) \in \textbf{M}_{\lambda})}\bigg\},
\label{eqn:G}
\end{eqnarray}
the LHS of (\ref{eqn:E=det}) is expanded as
\begin{equation}
\frac{1}{\prod_{m=1}^{M} N_m !} 
\sum_{\sigma \in \cS_N} (-1)^{N-\ell(\sigma)}
\prod_{\lambda=1}^{\ell(\sigma)}
\sum_{\substack{\textbf{M}_{\lambda}: \\
\{\sc_{\lambda}\} \setminus \cC(\sc_{\lambda}) 
\subset \textbf{M}_{\lambda} \subset \{\sc_{\lambda}\}}}
(-1)^{\sharp \textbf{M}^{\rm c}_{\lambda}}
G(\sc_{\lambda}, \textbf{M}_{\lambda}).
\label{eqn:LHS2}
\end{equation}
From now on, we will explain how to 
rewrite $G(\sc_{\lambda}, \textbf{M}_{\lambda})$ until (\ref{eqn:G4}).
We note that if we set
\begin{eqnarray}
&& F(\{x_{c_{\lambda}(j)}:c_{\lambda}(j) \in \textbf{M}_{\lambda}^{\rm c} \}) 
\nonumber\\
&& \qquad
= \int_{\mathbb{R}^{\textbf{M}_{\lambda}}}
\prod_{i: c_{\lambda}(i) \in \textbf{M}_{\lambda}}\Big\{ dx_{c_{\lambda}(i)} 
\ \chi_{{\tau}_{c_{\lambda}(i)}}(x_{c_{\lambda}(i)})
\sG_{{\tau}_{c_{\lambda}(i)}, {\tau}_{c_{\lambda}(i+1)}}
(x_{c_{\lambda}(i)},x_{c_{\lambda}(i+1)}) \Big\}
\nonumber\\
&& \qquad \qquad \times
\prod_{j: c_{\lambda}(j) \in \textbf{M}_{\lambda}^{\rm c}}
p_{{\tau}_{c_{\lambda}(j+1)}, {\tau}_{c_{\lambda}(j)}}
(x_{c_{\lambda}(j+1)},x_{c_{\lambda}(j)}),
\label{eqn:F}
\end{eqnarray}
which is the integral over $\mathbb{R}^{\textbf{M}_{\lambda}}$, 
then (\ref{eqn:G}) is obtained by performing the integral of it
over $\mathbb{R}^{\textbf{M}^{\rm c}_{\lambda}}=\mathbb{R}^{\{\sc_{\lambda}\}}
\setminus \mathbb{R}^{\textbf{M}_{\lambda}}$, 
\begin{equation}
G(\sc_{\lambda}, \textbf{M}_{\lambda}) =
\int_{\mathbb{R}^{\textbf{M}_{\lambda}^{\rm c}}} 
\prod_{j:c_{\lambda}(j) \in \textbf{M}_{\lambda}^{\rm c}} \Big\{ dx_{c_{\lambda}(j)} 
\chi_{{\tau}_{c_{\lambda}(j)}}(x_{c_{\lambda}(j)}) \Big\}
F(\{x_{c_{\lambda}(j)}:c_{\lambda}(j) \in \textbf{M}_{\lambda}^{\rm c} \}).
\label{eqn:G2}
\end{equation}
In (\ref{eqn:F}), use the integral representation 
(\ref{eqn:sGtt}) for
$\sG_{{\tau}_{c_{\lambda}(i)}, {\tau}_{c_{\lambda}(i+1)}}
(x_{c_{\lambda}(i)},x_{c_{\lambda}(i+1)})$
by putting the integral variables
to be $v=v_{c_{\lambda}(i)}$ and $w=w_{c_{\lambda}(i+1)}$. We obtain
\begin{eqnarray}
&& 
F(\{x_{c_{\lambda}(j)}:c_{\lambda}(j) \in \textbf{M}_{\lambda}^{\rm c}\}) 
\nonumber\\
&=& \int_{\mathbb{R}^{ \textbf{M}_{\lambda}}}
\prod_{i: c_{\lambda}(i) \in \textbf{M}_{\lambda}} \xi (dv_{c_{\lambda}(i)})
\int_{\mathbb{R}^{ \textbf{M}_{\lambda}}}
\prod_{i: c_{\lambda}(i) \in \textbf{M}_{\lambda}} \Big\{ dx_{c_{\lambda}(i)}
p_{0, {\tau}_{c_{\lambda}(i)}}(v_{c_{\lambda}(i)},x_{c_{\lambda}(i)}) 
\chi_{{\tau}_{c_{\lambda}(i)}}(x_{c_{\lambda}(i)}) \Big\}
\nonumber\\
&\times& \int_{\mathbb{R}^{ \textbf{M}_{\lambda}}}
\prod_{i:c_{\lambda}(i) \in \textbf{M}_{\lambda}} \Big\{ d w_{c_{\lambda}(i+1)}
p_{0, {\tau}_{c_{\lambda}(i+1)}}(0, w_{c_{\lambda}(i+1)})
\Phi_{\xi}^{v_{c_{\lambda}(i)}}(x_{c_{\lambda}(i+1)}+\sqrt{-1} w_{c_{\lambda}(i+1)}) \Big\}
\nonumber\\
&& \qquad \times
\prod_{j: c_{\lambda}(j) \in \textbf{M}_{\lambda}^{\rm c}} 
p_{{\tau}_{c_{\lambda}(j+1)}, {\tau}_{c_{\lambda}(j)}}
(x_{c_{\lambda}(j+1)}, x_{c_{\lambda}(j)}).
\nonumber\\
&=&
\int_{\mathbb{R}^{ \textbf{M}_{\lambda}}}
\prod_{i: c_{\lambda}(i) \in \textbf{M}_{\lambda}} \xi (dv_{c_{\lambda}(i)})
\textbf{E}_{\mib{v}} \Bigg[ 
\prod_{i:c_{\lambda}(i) \in \textbf{M}_{\lambda}} 
\bigg\{ \chi_{{\tau}_{c_{\lambda}(i)}}(V_{c_{\lambda}(i)}({\tau}_{c_{\lambda}(i)}))
\nonumber\\
&& \qquad \qquad \qquad \times
\Phi_{\xi}^{v_{c_{\lambda}(i)}}(Z_{c_{\lambda}(i+1)}
({\tau}_{c_{\lambda}(i+1)}))^{\1(c_{\lambda}(i+1) \in \textbf{M}_{\lambda})}
\nonumber\\
&& \qquad \qquad \qquad  \times 
\Phi_{\xi}^{v_{c_{\lambda}(i)}}
(x_{c_{\lambda}(i+1)}+\sqrt{-1}W_{c_{\lambda}(i+1)}
({\tau}_{c_{\lambda}(i+1)}))^{\1(c_{\lambda}(i+1) 
\in \textbf{M}_{\lambda}^{\rm c})}\bigg\}
\nonumber\\
&& \quad \times
\prod_{j: c_{\lambda}(j) \in \textbf{M}_{\lambda}^{\rm c}} 
\bigg\{ p_{{\tau}_{c_{\lambda}(j+1)}, {\tau}_{c_{\lambda}(j)}}
(V_{c_{\lambda}(j+1)}({\tau}_{c_{\lambda}(j+1)}), 
x_{c_{\lambda}(j)})^{\1(c_{\lambda}(j+1) \in \textbf{M}_{\lambda})} \nonumber\\
&& \qquad \qquad \qquad \qquad \times 
p_{{\tau}_{c_{\lambda}(j+1)}, {\tau}_{c_{\lambda}(j)}}
(x_{c_{\lambda}(j+1)}, x_{c_{\lambda}(j)})
^{\1(c_{\lambda}(j+1) \in \textbf{M}_{\lambda}^{\rm c})} \bigg\} \Bigg].
\nonumber
\end{eqnarray}
Using Fubini's theorem, (\ref{eqn:G2}) is given by
\begin{eqnarray}
&& \int_{\mathbb{R}^{\textbf{M}_{\lambda}}}
\prod_{i:c_{\lambda}(i) \in \textbf{M}_{\lambda}} \xi(dv_{c_{\lambda}(i)})
\textbf{E}_{\mib{v}} \Bigg[ 
\prod_{i:c_{\lambda}(i) \in \textbf{M}_{\lambda}} 
\chi_{{\tau}_{c_{\lambda}(i)}}(V_{c_{\lambda}(i)}({\tau}_{c_{\lambda}(i)}))
\nonumber\\
&& \qquad \qquad \times 
\prod_{i: c_{\lambda}(i), c_{\lambda}(i+1) \in \textbf{M}_{\lambda}}
\Phi_{\xi}^{v_{c_{\lambda}(i)}}(Z_{c_{\lambda}(i+1)}
({\tau}_{c_{\lambda}(i+1)}))
\nonumber\\
&& \quad \times 
\int_{\mathbb{R}^{\textbf{M}_{\lambda}^{\rm c}}} 
\prod_{j: c_{\lambda}(j) \in \textbf{M}_{\lambda}^{\rm c}} 
\Big\{ dx_{c_{\lambda}(j)} \chi_{\tau_{c_{\lambda}(j)}}
(x_{c_{\lambda}(j)}) \Big\}
\nonumber\\
&& \qquad \qquad \times 
\prod_{j: c_{\lambda}(j) \in 
\textbf{M}_{\lambda}^{\rm c}, c_{\lambda}(j+1) \in \textbf{M}_{\lambda}}
p_{{\tau}_{c_{\lambda}(j+1)}, {\tau}_{c_{\lambda}(j)}}
(V_{c_{\lambda}(j+1)}(\tau_{c_{\lambda}(j+1)}), x_{c_{\lambda}(j)})
\nonumber\\
&& \qquad \qquad \times
\prod_{j: c_{\lambda}(j), c_{\lambda}(j+1) \in \textbf{M}_{\lambda}^{\rm c}}
p_{{\tau}_{c_{\lambda}(j+1)}, {\tau}_{c_{\lambda}(j)}}
(x_{c_{\lambda}(j+1)}, x_{c_{\lambda}(j)})
\nonumber\\
&& \qquad \qquad \times
\prod_{i: c_{\lambda}(i) \in \textbf{M}_{\lambda}, 
c_{\lambda}(i+1) \in \textbf{M}_{\lambda}^{\rm c}}
\Phi_{\xi}^{v_{c_{\lambda}(i)}}
(x_{c_{\lambda}(i+1)}+\sqrt{-1} W_{c_{\lambda}(i+1)}
(\tau_{c_{\lambda}(i+1)}))
\Bigg].
\label{eqn:G3}
\end{eqnarray}
For each $1 \leq i \leq q_{\lambda}$ with 
$c_{\lambda}(i) \in \textbf{M}_{\lambda}$, we define
$\overline{i}=\min\{j > i : c_{\lambda}(j) \in \textbf{M}_{\lambda} \}$ 
and
$\underline{i}=\max\{j < i : c_{\lambda}(j) \in \textbf{M}_{\lambda} \}$. 
Then we perform integration over $x_{c_{\lambda}(j)}$'s
for $c_{\lambda}(j) \in \textbf{M}^{\rm c}_{\lambda}$
before taking the expectation $\textbf{E}_{\mib{v}}$.
That is, integrals over $x_{c_{\lambda}(j)}$'s
with indices in intervals $\underline{i} < j < i$
for all $i$, s.t. $c_{\lambda}(i) \in \textbf{M}_{\lambda}$
are done.
For each $i$, s.t. $c_{\lambda}(i) \in \textbf{M}_{\lambda}$,
if $\underline{i} < i-1$,
\begin{eqnarray}
&& \chi_{{\tau}_{c_{\lambda}(i)}}(V_{c_{\lambda}(i)}({\tau}_{c_{\lambda}(i)}))
 \prod_{j=\underline{i}+1}^{i-1}
\int_{\mathbb{R}} dx_{c_{\lambda}(j)} \chi_{\tau_{c_{\lambda}(j)}}
(x_{c_{\lambda}(j)}) 
p_{\tau_{c_{\lambda}(i)}, \tau_{c_{\lambda}(i-1)}}
(V_{c_{\lambda}(i)}(\tau_{c_{\lambda}(i)}), x_{c_{\lambda}(i-1)})
\nonumber\\
&& \qquad \quad \times
\prod_{k=\underline{i}+2}^{i-1} 
p_{\tau_{c_{\lambda}(k)}, \tau_{c_{\lambda}(k-1)}}
(x_{c_{\lambda}(k)}, x_{c_{\lambda}(k-1)})
\Phi_{\xi}^{v_{c_{\lambda}(\, \underline{i} \,)}}
(x_{c_{\lambda}(\underline{i}+1)}+\sqrt{-1} W_{c_{\lambda}(\underline{i}+1)}
(\tau_{c_{\lambda}(\underline{i}+1)})) 
\nonumber
\end{eqnarray}
coincides with the conditional expectation of
\begin{eqnarray}
&&\prod_{j=\underline{i}+1}^{i} \chi_{\tau_{c_{\lambda}(j)}}
(V_{c_{\lambda}(i)}(\tau_{c_{\lambda}(j)}))
\Phi_{\xi}^{v_{c_{\lambda}(\, \underline{i} \,)}}
(V_{c_{\lambda}(i)}(\tau_{c_{\lambda}(\underline{i}+1)})
+\sqrt{-1} W_{c_{\lambda}(\underline{i}+1)}
(\tau_{c_{\lambda}(\underline{i}+1)})).
\nonumber
\end{eqnarray}
with respect to $\textbf{E}_{\mib{v}}[\cdot|V_{c_\lambda(i)},W_{c_{\lambda}(\underline{i}+1)}]$.
Since $W_{i}(\cdot), i \in \{\sc_{\lambda} \}$ are i.i.d. random variables 
which are independent of $V_{i}(\cdot), i \in \{\sc_{\lambda} \}$,
$V_{c_{\lambda}(i)}(\tau_{c_{\lambda}(\underline{i}+1)})
+\sqrt{-1} W_{c_{\lambda}(\underline{i}+1)}$
$(\tau_{c_{\lambda}(\underline{i}+1)})$
has the same distribution as
$
V_{c_{\lambda}(i)}(\tau_{c_{\lambda}(\underline{i}+1)})
+\sqrt{-1} W_{c_{\lambda}(i)}
(\tau_{c_{\lambda}(\underline{i}+1)})
=Z_{c_{\lambda}(i)}(\tau_{c_{\lambda}(\underline{i}+1)}).
$ \\
Since 
$\prod_{i:c_{\lambda}(i) \in \textbf{M}_{\lambda}}
\Phi_{\xi}^{v_{c_{\lambda}(\, \underline{i} \,)}}$ 
$(Z_{c_{\lambda}(i)}(\tau_{c_{\lambda}(\underline{i}+1)}))
= \prod_{i:c_{\lambda}(i) \in \textbf{M}_{\lambda}}
\Phi_{\xi}^{v_{c_{\lambda}(i)}}
(Z_{c_{\lambda}(\, \overline{i} \,)}(\tau_{c_{\lambda}(i+1)}))$,
(\ref{eqn:G3}) is equal to 
\[
\int_{\mathbb{R}^{\textbf{M}_{\lambda}}} 
\prod_{i: c_{\lambda}(i) \in \textbf{M}_{\lambda}} \xi(dv_{c_{\lambda}(i)})
\textbf{E}_{\mib{v}} \left[ 
\prod_{i: c_{\lambda}(i) \in \textbf{M}_{\lambda}}
\prod_{j=\underline{i}+1}^{i}
\chi_{\tau_{c_{\lambda}(j)}}(V_{c_{\lambda}(i)}(\tau_{c_{\lambda}(j)})) 
\Phi_{\xi}^{v_{c_{\lambda}(i)}}
(Z_{c_{\lambda}(\, \overline{i} \,)}(\tau_{c_{\lambda}(i+1)}))
\right]. 
\]
Using only the entries of $\textbf{M}_{\lambda}$,
we can define a subcycle $\widehat{\sc_{\lambda}}$ of $\sc_{\lambda}$
uniquely as follows:
Since $\sc_{\lambda}$ is a cyclic permutation,
$\widehat{q_{\lambda}} \equiv \sharp \textbf{M}_{\lambda} \geq 1$.
Let $i_1=\min \{1 \leq i \leq q_{\lambda}: c_{\lambda}(i) 
\in \textbf{M}_{\lambda}\}$.
If $\widehat{q_{\lambda}} \geq 2$, define
$i_{j+1}=\overline{i}_j, 1 \leq j \leq \widehat{q_{\lambda}}-1$.
Then $\widehat{\sc_{\lambda}} 
=(\widehat{c_{\lambda}}(1) \widehat{c_{\lambda}}(2) 
\cdots \widehat{c_{\lambda}}(\widehat{q_{\lambda}}))
\equiv (c_{\lambda}(i_1) c_{\lambda}(i_2) 
\cdots c(i_{\widehat{q_{\lambda}}}))$.
Moreover, 
we decompose the set $\textbf{M}_{\lambda}$ into $M$ subsets,
$\textbf{M}_{\lambda}=\bigcup_{m=1}^{M} \textbf{J}_m^{\lambda}$,
by letting 
\[
\textbf{J}_m^{\lambda}=\textbf{J}_m^{\lambda}(\sc_{\lambda}, \textbf{M}_{\lambda})=
\Big\{ c_{\lambda}(i) \in \textbf{M}_{\lambda} :
\underline{i} < ^{\exists}j \leq i, \,
\mbox{s.t.} \, c_{\lambda}(j) \in \mathbb{I}^{(m)} \Big\}, \quad 1 \leq m \leq M.
\]
Note that by definition
$\textbf{J}_m^{\lambda} \cap \textbf{J}_{m'}^{\lambda} \not= \emptyset,
m \not= m'$ in general, and
$\textbf{J}_1^{\lambda}=\mathbb{I}_{N_1} \cap \textbf{M}_{\lambda}=
\mathbb{I}_{N_1} \cap \{\sc_{\lambda}\}$, 
$\textbf{J}_m^{\lambda} \subset \mathbb{I}_{\sum_{k=1}^{m} N_k}$ for $2 \leq m \leq M$,
$\textbf{J}_m^{\lambda} \cap \mathbb{I}^{(k)} \subset \textbf{J}_{k}^{\lambda}$
for $1 \leq k < m \leq M$.
Finally we arrive at the following expression
of $G(\sc_{\lambda}, \textbf{M}_{\lambda})$,
\begin{equation}
\int_{\mathbb{R}^{\textbf{M}_{\lambda}}} \prod_{i:c_{\lambda}(i) \in \textbf{M}_{\lambda}}
\xi(dv_{c_{\lambda}(i)}) \
\textbf{E}_{\mib{v}} \left[ 
\prod_{m=1}^{M} \prod_{j_m \in \textbf{J}_m^{\lambda}}
\chi_{t_m}(V_{j_m}(t_m))
\prod_{i=1}^{\widehat{q_{\lambda}}}
\Phi_\xi^{v_{\widehat{c_{\lambda}}(i)}} 
(Z_{\widehat{c_{\lambda}}(i+1)}(T)) \right],
\label{eqn:G4}
\end{equation}
where the martingale property (\ref{eqn:martingale})
is used.
Let $\textbf{M} \equiv \bigcup_{\lambda=1}^{\ell(\sigma)}
\textbf{M}_{\lambda}$. 
Since
$N-\sum_{\lambda=1}^{\ell(\sigma)} \sharp \textbf{M}_{\lambda}^{\rm c}
=\sharp \textbf{M}$,
the LHS of (\ref{eqn:E=det}), which is written as 
(\ref{eqn:LHS2}), becomes now
\begin{eqnarray}
&& \frac{1}{\prod_{m=1}^{M} N_m !} 
\sum_{\sigma \in \cS_N} 
\sum_{\substack{\textbf{M} : \\
\mathbb{I}_N \setminus \bigcup_{\lambda=1}^{\ell(\sigma)}
\cC(\sc_{\lambda}) \subset \textbf{M} \subset \mathbb{I}_N}}
(-1)^{\sharp \textbf{M} -\ell(\sigma)}
\int_{\mathbb{R}^{\textbf{M}}}
\prod_{\lambda=1}^{\ell(\sigma)}
\prod_{i: c_{\lambda}(i) \in \textbf{M}_{\lambda}}
\xi(dv_{c_{\lambda}(i)})
\nonumber\\
&& \qquad \qquad \times 
\textbf{E}_{\mib{v}} \left[ \prod_{\lambda=1}^{\ell(\sigma)}
\left\{ \prod_{m=1}^{M}
\prod_{j_m \in \textbf{J}_m^{\lambda}} 
\chi_{t_m}(V_{j_m}(t_m)) 
\prod_{i=1}^{\widehat{q_{\lambda}}}
\Phi_{\xi}^{v_{\widehat{c_{\lambda}}(i)}}
(Z_{\widehat{c_{\lambda}}(i+1)}(T)) 
\right\} \right].
\label{eqn:E=det2}
\end{eqnarray}

We define $\widehat{\sigma} \equiv 
\widehat{\sc_1} \widehat{\sc_2} \cdots
\widehat{\sc_{\ell(\sigma)}}$
and  $\textbf{J}_m \equiv \bigcup_{\lambda=1}^{\ell(\sigma)} \textbf{J}_m^{\lambda},
1 \leq m \leq M$.
Note that $\ell(\widehat{\sigma})=\ell(\sigma)$.
The obtained $(\textbf{J}_m)_{m=1}^{M}$'s 
form a collection of series of index sets
satisfying the following conditions,
which we write as ${\cal J}(\{N_m\}_{m=1}^{M})$:

%%%%%%%%%%%%%%%%%%%%%%%%%%%%%%%%%%%%%%%%%%%%%%%%%%%%%%
\vskip 0.3cm
\noindent (\textbf{C.J}) \,
$\textbf{J}_1=\mathbb{I}_{N_1}, \, 
\textbf{J}_m \subset \mathbb{I}_{\sum_{k=1}^{m} N_{k}}$ 
for $2 \leq m \leq M$, \, 
$\textbf{J}_m \cap \mathbb{I}^{(k)} \subset \textbf{J}_{k}$ for
$1 \leq k < m \leq M$, \mbox{and} \\
\qquad \quad $\sharp \textbf{J}_m=N_m$ for $1 \leq m \leq M$.
\vskip 0.3cm
%%%%%%%%%%%%%%%%%%%%%%%%%%%%%%%%%%%%%%%%%%%%%%%%%%%%%%

For each $(\textbf{J}_m)_{m=1}^{M} \in {\cal J}(\{N_m\}_{m=1}^{M})$,
we put $A_1=0$ and
$A_m=\sharp(\textbf{J}_m \cap \mathbb{I}_{\sum_{k=1}^{m-1}N_{k}})
=\sharp ( \textbf{J}_m \cap \bigcup_{k=1}^{m-1} \textbf{J}_k), 2 \leq m \leq M$.
Then, if we put $\textbf{M}=\bigcup_{m=1}^{M} \textbf{J}_m$,
$\sharp \textbf{M}=\sum_{m=1}^{M}(N_m-A_m)$,
which means that
from the original index set $\mathbb{I}_N=\bigcup_{m=1}^{M} \mathbb{I}^{(m)}$
with $\sharp \mathbb{I}^{(m)}=N_m, 1 \leq m \leq M$,
we obtain a subset $\textbf{M}$ by eliminating
$A_m$ elements at each level $1 \leq m \leq M$.
By this reduction, we obtain
$\widehat{\sigma} \in \cS(\textbf{M})$ from $\sigma \in \cS_N$.
It implies that, for all $\widehat{\sigma} \in \cS(\textbf{M})$,
the number of $\sigma$'s in $\cS_N$ which give
the same $\widehat{\sigma}$ and
$(\textbf{J}_m)_{m=1}^{M}$ by this reduction is given by
$\prod_{m=1}^{M} A_m !$, where $0! \equiv 1$.
Then (\ref{eqn:E=det2}) is equal to 
\begin{eqnarray}
&& 
\sum_{\substack{\textbf{M}: \\
\max_m\{N_m\} \leq \sharp \textbf{M} \leq N}}
\sum_{\substack{(\textbf{J}_m)_{m=1}^{M} \subset
{\cal J}(\{N_m\}_{m=1}^{M}): \\
\bigcup_{m=1}^{M} \textbf{J}_m=\textbf{M}
}}
\frac{\prod_{m=1}^{M} A_m !}{\prod_{m=1}^{M} N_m !}
\sum_{\widehat{\sigma} \in \cS(\textbf{M})}
(-1)^{\sharp \textbf{M}-\ell(\widehat{\sigma})}
\nonumber\\
&& \qquad \times
\sharp \textbf{M} !
\int_{\mathbb{W}^{\rm A}_{\sharp \textbf{M}}} \xi^{\otimes \textbf{M}} (d \mib{v})
\textbf{E}_{\mib{v}} \left[
\prod_{m=1}^{M} \prod_{j_m \in \textbf{J}_m}
\chi_{t_m}(V_{j_m}(t_m))
\prod_{\lambda=1}^{\ell(\widehat{\sigma})}
\prod_{i=1}^{\widehat{q_{\lambda}}}
\Phi_{\xi}^{v_{\widehat{ c_{\lambda} }(i)}}
(Z_{\widehat{c_{\lambda}}(i+1)}(T)) \right]
\nonumber\\
&&=
\sum_{\substack{\textbf{M}: \\
\max_m\{N_m\} \leq \sharp \textbf{M} \leq N}}
\sum_{\substack{(\textbf{J}_m)_{m=1}^{M} \subset
{\cal J}(\{N_m\}_{m=1}^{M}): \\
\bigcup_{m=1}^{M} \textbf{J}_m=\textbf{M}
}}
\sharp \textbf{M} ! \prod_{m=1}^{M} \frac{A_m!}{N_m!}
\nonumber\\
&& \qquad \times
\int_{\mathbb{W}^{\rm A}_{\sharp \textbf{M}}} \xi^{\otimes \textbf{M}} (d \mib{v})
\textbf{E}_{\mib{v}} \left[
\prod_{m=1}^{M} \prod_{j_m \in \textbf{J}_m}
\chi_{t_m}(V_{j_m}(t_m))
\det_{i,j \in \textbf{M}}
\left[ \Phi_{\xi}^{v_i} (Z_{j}(T)) \right] \right].
\label{eqn:AAA}
\end{eqnarray}
Assume $\max_{m}\{N_m\} \leq p \leq N$, 
$0 \leq A_m \leq N_m, 2 \leq m \leq M$ and set $A_1 =0$.
Consider
\begin{eqnarray}
&& \Lambda_1 =\left\{(\textbf{J}_m)_{m=1}^{M} \subset
{\cal J}(\{N_{m}\}_{m=1}^{M}): 
\sharp (\bigcup_{m=1}^{M} \textbf{J}_m) =p,
\right. 
\nonumber\\
&& \hskip 6cm \left.
\sharp (\textbf{J}_m \cap \bigcup_{k=1}^{m-1} \textbf{J}_k)
=A_m, 2 \leq m \leq M \right\},
\nonumber\\
&& \Lambda_2 = \left\{(\mathbb{J}_m)_{m=1}^{M} :
\sharp \mathbb{J}_m=N_m, 1 \leq m \leq M,
\bigcup_{m=1}^{M} \mathbb{J}_m= \mathbb{I}_{p},
\right. 
\nonumber\\
&& \hskip 6cm \left.
\sharp (\mathbb{J}_m \cap \bigcup_{k=1}^{m-1} \mathbb{J}_k )
=A_m, 2 \leq m \leq M \right\}.
\nonumber
\end{eqnarray}
Since the CBMs are i.i.d. in $\textbf{P}_{\mib{v}}$, 
the integral in (\ref{eqn:AAA}) has the same value
for all $(\textbf{J}_m)_{m=1}^{M} \in \Lambda_1$
with $\bigcup_{m=1}^{M} \textbf{J}_{m}=\textbf{M}$
and it is also equal to
\[
\int_{\mathbb{W}^{\rm A}_{p}} \xi^{\otimes p} (d \mib{v})
\textbf{E}_{\mib{v}} \left[
\prod_{m=1}^{M} \prod_{j_m \in \mathbb{J}_m}
\chi_{t_m}(V_{j_m}(t_m))
\det_{i,j \in \mathbb{I}_{p}}
\left[ \Phi_{\xi}^{v_i} (Z_{j}(T)) \right] \right]
\]
for all $(\mathbb{J}_{m})_{m=1}^{M} \in \Lambda_2$.
In $\Lambda_1$, for each $2 \leq m \leq M$,
$A_m$ elements in $\textbf{J}_m$ are chosen from
$\bigcup_{k=1}^{m-1} \textbf{J}_k$,
in which $\sharp (\bigcup_{k=1}^{m-1}
\textbf{J}_{k})=\sum_{k=1}^{m-1}(N_k-A_k)$,
and the remaining $N_m-A_m$ elements in $\textbf{J}_m$
are from $\mathbb{I}^{(m)}$ with $\sharp \mathbb{I}^{(m)}=N_m$.
Then
\[
\sharp \Lambda_1=\prod_{m=2}^{M} 
\binom{\sum_{k=1}^{m-1} (N_k-A_k)}{A_m}
\binom{N_m}{N_m-A_m}.
\]
In $\Lambda_2$, on the other hand,
$N_1$ elements in $\mathbb{J}_1$ is chosen from $\mathbb{I}_{p}$,
and then for each $2 \leq m \leq M$,
$A_m$ elements in $\mathbb{J}_m$ are chosen from
$\bigcup_{k=1}^{m-1} \mathbb{J}_k$ with
$\sharp (\bigcup_{k=1}^{m-1} \mathbb{J}_k)=
\sum_{k=1}^{m-1}(N_k-A_k)$
and the remaining $N_m-A_m$ elements in $\mathbb{J}_m$
are from $\mathbb{I}_{p} \setminus \bigcup_{k=1}^{m-1} \mathbb{J}_k$
with $\sharp(\mathbb{I}_{p} \setminus \bigcup_{k=1}^{m-1} \mathbb{J}_k)
=p-\sum_{k=1}^{m-1}(N_k-A_k)$.
Then
\[
\sharp \Lambda_2
=\binom{p}{N_1} 
\prod_{m=2}^{M} 
\binom{\sum_{k=1}^{m-1} (N_k-A_k)}{ A_m}
\binom{p-\sum_{k=1}^{m-1} (N_k-A_k)}{ N_m-A_m}.
\]
Since $\sum_{m=1}^{M}(N_m-A_m)=p$,
we see 
$\sharp \Lambda_2/\sharp \Lambda_1=
p ! \prod_{m=1}^{M} A_m !/N_m!$.
Then (\ref{eqn:AAA}) is equal to the RHS of (\ref{eqn:E=det})
and the proof is completed.
\hfill $\Box$

%%%%%%%%%%%%%%%%%%%%%%%%%%%%%%%%%%%%%%%%%%%%%%%%%%%%%%%%%%
%%%%%%%%%%%%%%%%%%%%%%%%%%%%%%%%%%%%%%%%%%%%%%%%%%%%%%%%%%
%%%%%%%%%%%%%%%%%%%%%%%%%%%%%%%%%%%%%%%%%%%%%%%%%%%%%%%%%%
%%%  SEC4   %%%%%%%%%%%%%%%%%%%%%%%%%%%%
%%%%%%%%%%%%%%%%%%%%%%%%%%%%%%%%%%%%%%%%%%%%%%%%%%%%%%%%%%
\SSC{Proof of Corollary \ref{cor:infinite_cBM}}%%%%%%%%%
%%%%%%%%%%%%%%%%%%%%%%%%%%%%%%%%%%%%%%%%%%%%%%%%%%%%%%%%%%

We consider the case that $k=1$ and $G(x)=x^q$,
$q \in \mathbb{N}$.
(The argument will be easily extended to general polynomials
of order $q$.)
We introduce a map $\pi$ form $\bigoplus_{n=1}^\infty \mathbb{R}^n$ 
to $\bigoplus_{p=1}^\infty \mathbb{W}_{p}^{\rm A}$
such that
\[
\pi (w_1,w_2,\dots,w_n) = (v_1, v_2, \dots, v_p),
\quad \mbox{where $\{w_i\}_{i=1}^n = \{v_i\}_{i=1}^p$}.
\]
We also introduce the functions
$
p(\mib{w})=\sharp \{w_i\}_{i=1}^n
$, and
$a_i(\mib{w})= \sharp \{j : w_j= (\pi \mib{w})_i \}, i \in \mathbb{I}_{p(\mib{w})}$. 
Then in case $\xi(\mathbb{R})<\infty$ 
the CBM representation (\ref{eqn:main}) gives
\begin{eqnarray}
&&\mathbb{E}_{\xi}\left[F\big(\Xi(\cdot)\big)\right]
=\textbf{E}_{\mib{u}} \left[\left( \sum_{i=1}^{\xi(\mathbb{R})} \phi_1(V_i(t_1)) \right)^q
\det_{1\le i,j \le {\xi(\mathbb{R})}}
\Big[\Phi_{\xi}^{u_i}(Z_j(T))\Big]\right]
\nonumber\\
&&=\int_{\mathbb{R}^q}\xi^{\otimes q}(d\mib{w})
\textbf{E}_{\pi \mib{w}} \left[\prod_{i=1}^{p(\mib{w})} \phi_1(V_i(t_1))^{a_i(\mib{w})}
\det_{i,j \in \mathbb{I}_{p (\mib{w})}}
\Big[\Phi_{\xi}^{(\pi \mib{w})_i}(Z_j(T))\Big]\right]
\nonumber\\
&&= \sum_{p=1}^q 
\sum_{\substack{a_i\in\mathbb{N}, i \in \mathbb{I}_{p}: \\
\sum_{i=1}^p a_i=q}}
\binom{q}{ a_1,a_2,\dots,a_p}
\int_{\mathbb{W}^{\rm A}_{p}}\xi^{\otimes p}(d\mib{v})
\textbf{E}_{\mib{v}} \left[\prod_{i=1}^{p} \phi_1(V_i(t_1))^{a_i}
\det_{i,j \in \mathbb{I}_{p}}
\Big[\Phi_{\xi}^{v_i}(Z_j(T))\Big]\right].
\nonumber\\
\label{eqn:mainB}
\end{eqnarray}
Here we used the fact that for $\mib{v}\in \mathbb{W}^{\rm A}_{p}$,
$a_i \in \mathbb{N}, 1 \leq i \leq p, \sum_{i=1}^{p} a_i=q$, 
\[
\sharp \Big\{\mib{w} \in \mathbb{R}^q: \pi \mib{w}=\mib{v}, \ 
a_i(\mib{w})=a_i, 1\le i \le p \Big\}
=\frac{q!}{a_1! a_2! \cdots a_p !}
\equiv \binom{q}{ a_1,a_2,\dots,a_p},
\]
and the equality in the measure $\xi^{\otimes p}(d \mib{v})$
\begin{eqnarray}
&&\textbf{E}_{\mib{v}} \left[ \prod_{i=1}^{p} \phi_1(V_i(t_1))^{a_i}
\det_{1\le i,j \le {\xi(\mathbb{R})}}
\Big[\Phi_{\xi}^{v_i}(Z_j(T))\Big]\right]
\nonumber\\
&&\qquad\qquad=\textbf{E}_{\mib{v}} \left[ \prod_{i=1}^{p} \phi_1(V_i(t_1))^{a_i}
\det_{i,j \in \mathbb{I}_p}
\Big[\Phi_{\xi}^{v_i}(Z_j(T))\Big]\right],
\label{eqn:point1}
\end{eqnarray}
which holds for any $p \leq \xi(\mathbb{R})$ by (\ref{eqn:Kronecker2}).
Note that the equality (\ref{eqn:point1}) is valid
also in the case that $\xi(\mathbb{R})=\infty$ 
under the conditions $(\textbf{C.1})$ and $(\textbf{C.2})$.

By the bound (\ref{estimate:Phi}) obtained from 
the conditions $(\textbf{C.1})$ and $(\textbf{C.2})$,
we can prove the uniform integrability of the functions
$\displaystyle{
\prod_{i=1}^{p} \phi_1(V_i(t_1))^{a_i}
\1(|v_i|\le L)
\det_{i,j \in \mathbb{I}_{p}}
\Big[\Phi_{\xi}^{v_i}(Z_j(T))\Big]
}$, $L \in \mathbb{N}$, 
with respect to the measure $\xi^{\otimes p}(d\mib{v})\textbf{P}_{\mib{v}}$.
Then from (\ref{eqn:mainB}) we conclude that
\begin{eqnarray}
&& \lim_{L\to\infty}\mathbb{E}_{\xi\cap [-L,L]}\left[F\big(\Xi(\cdot)\big)\right]
%\nonumber\\
%&&
= \sum_{p=1}^q 
\sum_{\substack{a_i\in\mathbb{N}, i \in \mathbb{I}_{p}: \\
\sum_{i=1}^p a_i=q
}}
\binom{q}{ a_1,a_2,\dots,a_p}
\nonumber\\
&& \qquad \qquad \qquad \times
\int_{\mathbb{W}^{\rm A}_{p}}\xi^{\otimes p}(d\mib{v})
\textbf{E}_{\mib{v}} \left[\prod_{i=1}^{p} \phi_1(V_i(t_1))^{a_i}
\det_{i,j \in \mathbb{I}_{p}}
\Big[\Phi_{\xi}^{v_i}(Z_j(T))\Big]\right].
\label{eqn:expression1}
\end{eqnarray}
This is a realization of the RHS of (\ref{eqn:main}), 
when $F(\sum_{i=1}^{\xi(\mathbb{R})} \delta_{V_i(\cdot)})
=(\sum_{i=1}^{\xi(\mathbb{R})} \phi_1(V_{i}(t_1)))^q$.
If $q$ is finite, $q \in \mathbb{N}$, 
the sizes of matrices for the determinants
in (\ref{eqn:main}) can be
reduced from $\xi(\mathbb{R})$ to $p$ with $1 \leq p \leq q$
as in (\ref{eqn:expression1}). 
Then, even if $\xi(\mathbb{R})=\infty$, we needn't deal with 
infinite-dimensional determinants.
Generalization for $k \geq 2$ is straightforward. 
\hfill $\Box$

%%%%%%%%%%%%%%%%%%%%%%%%%%%%%%%%%%%%%%%%%%%%%%%%%%%%%%%%%%
%%%%%%%%%%%%%%%%%%%%%%%%%%%%%%%%%%%%%%%%%%%%%%%%%%%%%%%%%%
%%%%%%%%%%%%%%%%%%%%%%%%%%%%%%%%%%%%%%%%%%%%%%%%%%%%%%%%%%
%%%  SEC5   %%%%%\label{eqn:KolCri}%%%%%%%%%%%%%%%%%%%%%%%
%%%%%%%%%%%%%%%%%%%%%%%%%%%%%%%%%%%%%%%%%%%%%%%%%%%%%%%%%%
\SSC{Proof of Proposition \ref{prop:4-th_moment}}%%%%%%%%%
%%%%%%%%%%%%%%%%%%%%%%%%%%%%%%%%%%%%%%%%%%%%%%%%%%%%%%%%%%

Suppose that the initial configuration
$\xi\in \mM_0$ satisfies the conditions 
$(\textbf{C.1})$ and $(\textbf{C.2})$
with constants $C_0, C_1, C_2$ and indices $\alpha, \beta$.
From the proof of Corollary \ref{cor:infinite_cBM}
given in the previous section, 
we see that the LHS of (\ref{eqn:KolCri}) is given by 
\[
\sum_{p=1}^4 \int_{\mathbb{W}_p^{\rm A}}\xi^{\otimes p}(d\mib{v})
\textbf{E}_{\mib{v}}\Bigg[
F_{p}\bigg( \Big(\varphi(V_i(t))-\varphi(V_i(s)) 
\Big)_{i \in \mathbb{I}_{p}} \bigg)
\det_{1\le i,j \le p}\Big[\Phi_{\xi}^{v_i}(Z_j(T))\Big]
\Bigg]
\]
with
$F_1(x_1)= x_1^4$, 
$F_2(\mib{x}_2)= 6 x_1^2x_2^2+4 x_1 x_2 (x_1^2+x_2^2)$,
$F_3(\mib{x}_3)= 12x_1 x_2 x_3 (x_1+x_2+x_3)$
and
$F_4(\mib{x}_4)= 24 x_1 x_2 x_3 x_4$.
Then Proposition \ref{prop:4-th_moment} is 
concluded from the following estimate.

%%%%%%%%%%%%%%%%%%%%%%%%%%%%%%%%%%%%%%%%%%%%%%%%%%%%%%%%%%%
%%%%%%%%%%%%%%%%%%%%%%%%%%%%%%%%%%%%%%%%%%%%%%%%%%%%%%%%%%%
%% Lemma 5.1%%%%%%%%%%%%%%%%%%%%%%%%%%%%%%%%%%%%%%
%%%%%%%%%%%%%%%%%%%%%%%%%%%%%%%%%%%%%%%%%%%%%%%%%%%%%%%%%%%
\begin{lemma}
Let $\{a_i\}_{i=1}^p$ be a sequence of 
positive integers with length $p\in\mathbb{N}$.
Then for any $T>0$ and $\varphi \in C_0^{\infty}(\mathbb{R})$
there exists a positive constant 
$C=C(C_0,C_1,C_2,\alpha,\beta,T,\varphi)$,
which is independent of $s, t$, such that
\begin{eqnarray}
&&\int_{\mathbb{W}_p^{\rm A}}\xi^{\otimes p}(d\mib{v})
\textbf{E}_{\mib{v}}\left[
\prod_{i=1}^p \bigg(
\varphi(V_i(t))-\varphi(V_i(s)) \bigg)^{a_i}
\det_{1\le i,j \le p} \Big[\Phi_{\xi}^{v_i}(Z_j(T)) \Big]
\right]
\nonumber\\
&& \hskip 4cm
\le C |t-s|^{\sum_{i=1}^p a_i/2}, 
 \quad
 \forall s,t \in [0,T].
\label{est:moment}
\end{eqnarray}
\end{lemma}
%%%%%%%%%%%%%%%%%%%%%%%%%%%%%%%%%%%%%%%%%%%%%%%%%%%%%%%%%%%
%% Lemma %%%%%%%%%%%%%%%%%%%%%%%%%%%%%%%%%%%%%%%%
%%%%%%%%%%%%%%%%%%%%%%%%%%%%%%%%%%%%%%%%%%%%%%%%%%%%%%%%%%%

\noindent{\it Proof.} 
Choose $L \in \mathbb{N}$ so that $\textrm{supp } \varphi \subset [-L,L]$,
and put $\1_L(x,y)=0$, if $|x|>L$ and $|y|>L$, 
and $\1_L(x,y)$=1, otherwise.
By the Schwartz inequality the LHS of (\ref{est:moment}) 
is bounded from the above by
\begin{eqnarray}
&&\int_{\mathbb{W}_p^{\rm A}}\xi^{\otimes p}(d\mib{v})
\textbf{E}_{\mib{v}} \left[
\prod_{i=1}^p \bigg(
\varphi(V_i(t))-\varphi(V_i(s))\bigg)^{2a_i} \right]^{1/2}
\textbf{E}_{\mib{v}} \left[
\prod_{i=1}^p \1_L(V_i(s),V_i(t)) \right]^{1/4}
\nonumber\\
&&\qquad \qquad \quad \times \textbf{E}_{\mib{v}} \left[
\left(\det_{1\le i,j \le p} \Big[\Phi_{\xi}^{v_i}(Z_j(T))
\Big] \right)^4 \right]^{1/4}. 
\nonumber
\end{eqnarray}
Since $V_i(t)$, $i\in\mathbb{N}$ are independent Brownian motions, 
$
\textbf{E}_{\mib{v}} [
\prod_{i=1}^p (
\varphi(V_i(t))-\varphi(V_i(s)) )^{2a_i} ]
\le c_1|t-s|^{\sum_{i=1}^p a_i}
$
and
$
\textbf{E}_{\mib{v}} [
\prod_{i=1}^p \1_L(V_i(s),V_i(t)) ]
\le c_2 e^{ -c_2' \sum_{i=1}^p |v_i|^2 }, 
\forall s,t \in [0,T],
$
with positive constants $c_1, c_2, c_2'$, 
which are independent of $s,t$.
And from the estimate (\ref{estimate:Phi})
we have 
\[
\textbf{E}_{\mib{v}} \left[
\left(\det_{1\le i,j \le p} \Big[\Phi_{\xi}^{v_i}(Z_j(T))
\Big] \right)^4
\right]^{1/4}
\le c_3 \exp \bigg\{ c_3' \sum_{i=1}^p |v_i|^\theta \bigg\}
\]
with positive constants $c_3, c_3'$
and $\theta \in (\max\{\alpha, (2-\beta)\}, 2)$.
Combining the above estimates with the fact that, for any $c,c'>0$,
$
\int_{\mathbb{R}}\xi(dv) e^{c|v|^\theta -c'|v|^2} < \infty,
$
which is derived from the condition (\textbf{C.2}) (i)
and the fact $\theta<2$,
we obtain the lemma.
\hfill $\Box$

%%%%%%%%%%%%%%%%%%%%%%%%%%%%%%%%%%%%%%%%%%%%%%%%%%%%%%%%%%%
%%%%%%%%%%%%%%%%%%%%%%%%%%%%%%%%%%%%%%%%%%%%%%%%%%%%%%%%%%%
%%%%%%%%%%%%%%%%%%%%%%%%%%%%%%%%%%%%%%%%%%%%%%%%%%%%%%%%%%%
%%%  SEC6   %%%%%\label{eqn:KolCri}%%%%%%%%%%%%%%%%%%%%%%%%
%%%%%%%%%%%%%%%%%%%%%%%%%%%%%%%%%%%%%%%%%%%%%%%%%%%%%%%%%%%
\SSC{Proof of Proposition \ref{prop:noncolliding}}%%%%%%%%%
%%%%%%%%%%%%%%%%%%%%%%%%%%%%%%%%%%%%%%%%%%%%%%%%%%%%%%%%%%%

Let $\tau= \inf \{ t>0 : \Xi(t) \notin \mM_0 \}$
and $\tau_{ij}$ be defined by (\ref{eqn:tauij}).
From the CBM representation (\ref{eqn:main}), 
for any $\xi\in \mM_0$ with $\xi(\mathbb{R})\in\mathbb{N}$,
\begin{eqnarray}
\mathbb{P}_{\xi}\Big[\tau \le T \Big]
&\le& \textbf{E}_{\mib{u}} \left[ 
\sum_{1\le i < j \le \xi(\mathbb{R})}\1 (\tau_{ij}\le T) 
\det_{1\leq i,j \le \xi(\mathbb{R})}
\Big[\Phi_{\xi}^{u_i}(Z_j(T))\Big] 
\right]
\nonumber\\
&=& \int_{\mathbb{W}_2^{\rm A}} \xi^{\otimes 2}(d\mib{v}) \
\textbf{E}_{\mib{v}} \left[ \1 (\tau_{12}\le T) 
\det_{1\le i,j \le 2}
\Big[\Phi_{\xi}^{v_i}(Z_j(T))\Big] 
\right].
\label{ineq:tau}
\end{eqnarray}
Since $\theta$ in (\ref{estimate:Phi}) is strictly less than $2$, 
and the constants and the indices 
in the conditions for the configurations $\xi \cap [-L,L]$ 
can be taken to be independent of $L>0$, 
\[
\lim_{L\to\infty}\Phi_{\xi\cap [-L,L]}^{a}(Z_1(t)) = \Phi_{\xi}^{a}(Z_1(t))
\quad \mbox{ in $L^k(\Omega, \textbf{P}_{v})$}
\]
holds for any $k\in\mathbb{N}$.
Hence, the inequality (\ref{ineq:tau}) holds for 
$\xi\in \mM_0$ under the conditions $(\textbf{C.1})$ and $(\textbf{C.2})$.
By the strong Markov property of CBM 
started at $\mib{v}\in \mathbb{W}^{\rm A}_2$
\begin{eqnarray}
&&\textbf{E}_{\mib{v}} \left[ 
\1 (\tau_{12}\le T) 
\det_{1\le i,j \le 2}
\Big[\Phi_{\xi}^{v_i}(Z_j(T))\Big] 
\right]
\nonumber\\
&&\qquad\qquad
=\textbf{E}_{\mib{v}} \left[ 
\1 (\tau_{12}\le T) 
\textbf{E}_{\mib{Z}(\tau_{12})} \left[
\det_{1\le i,j \le 2}
\left[\Phi_{\xi}^{v_i}(Z_j(T-\tau_{12})) \right] \right]
\right].
\nonumber
\end{eqnarray}
By the martingale property of $\Phi_{\xi}^{v_i}(Z_j(T))$
we can apply the optional stopping theorem 
and show that the RHS of the above equation coincides with
\begin{eqnarray}
&&\textbf{E}_{\mib{v}} \left[ \1 (\tau_{12}\le T) 
\det_{1\le i,j \le 2} \left[
\Phi_{\xi}^{v_i}(Z_j(\tau_{12})) \right]
\right]
={\rm E}_1 \left[ \1 (\tau_{12}\le T)
{\rm E}_2 \left[\det_{1\le i,j \le 2}
\left[\Phi_{\xi}^{v_i}(Z_j(\tau_{12})) \right]
\right]\right]
\nonumber\\
&& \qquad =\frac{\sqrt{-1}}{v_1-v_2}
\prod_{r \in\textrm{supp }{\xi}\setminus\{v_1,v_2\}}
\frac{1}{(r-v_1)(r-v_2)}
\nonumber\\
&& \quad \qquad\times {\rm E}_1 \left[ \1 (\tau_{12}\le T)
{\rm E}_2 \left[ \big(W_1(\tau_{12})
-W_2(\tau_{12})\big)
\prod_{k=1,2}
G\big(V_k(\tau_{12}), W_k(\tau_{12})\big)
\right]\right],
\nonumber
\end{eqnarray}
where
$G(v,w)= \prod_{r\in \textrm{supp } \xi\setminus \{v_1,v_2\}} 
(r-v-\sqrt{-1}w)$,
and the fact that  $V_1(\tau_{12}) = V_2(\tau_{12})$ 
almost surely
was used in the last equality.
Since $W_{k}(\tau_{12}), k=1,2$ are i.i.d. under ${\rm P}_2$, we have
${\rm E}_2 \left[ \big(W_1(\tau_{12}) 
-W_2(\tau_{12})\big)
\prod_{k=1,2}
G\big(V_k(\tau_{12}),W_k(\tau_{12})\big)
\right]=0$.
This completes the proof.
\hfill $\Box$

%%%%%%%%%%%%%%%%%%%%%%%%%%%%%%%%%%%%%%%%%%%%%%%%%%%%%%%%%%%%%%%%%%%%
%%%%%%%%%%%%%%%%%%%%%%%%%%%%%%%%%%%%%%%%%%%%%%%%%%%%%%%%%%%%%%%%%%%%
\vskip 3mm
\begin{small}
%%%%%%%%%%%%%%%%%%%%%%%%%%%%%%%%%%%%%%%%%%%%%%%%%%%%%%
\noindent{\it Acknowledgements.} 
The present authors are grateful to the referee
for careful reading of the manuscript and useful comments.
A part of the present work was done
during the participation of the authors
in the international conference
``Selfsimilar Processes and Their Applications",
Angers, July 20-24, 2009,
whose last day program was prepared 
for the 60th birthday of Professor Marc Yor.
The authors would like to dedicate the present paper
to Professor Marc Yor.
M.K. is supported in part by
the Grant-in-Aid for Scientific Research (C)
(No.21540397) of Japan Society for
the Promotion of Science.
H.T. is supported in part by
the Grant-in-Aid for Scientific Research 
(KIBAN-C, No.19540114) of Japan Society for
the Promotion of Science.

%%%%%%%%%%%%%%%%%%%%%%%%%%%%%%%%%%%%%%%%%%%%%%%%%%%%%%
\end{small}
%%%%%%%%%%%%%%%%%%%%%%%%%%%%%%%%%%%%%%%%%%%%%%%%%%%%%%%%%%%%
%%%%% REFERENCES %%%%%%%%%%%%%%%%%%%%%%%%%%%%%%%%%%%%%%%%%%%
%%%%%%%%%%%%%%%%%%%%%%%%%%%%%%%%%%%%%%%%%%%%%%%%%%%%%%%%%%%%

\footnotesize 
%%%%%%%%%%%%%%%%%%%%%%%%%%%%%%%%%%%%%%%%%%%%%%%%%%%%%%%%%%%%%
%%%%%%%%%%%%%%%Reference%%%%%%%%%%%%%%%%%%%%%%%%%%%%%%%%%%%%%
%%%%%%%%%%%%%%%%%%%%%%%%%%%%%%%%%%%%%%%%%%%%%%%%%%%%%%%%%%%%%

%%%%%%%%%%%%%%%%%%%%%%%%%%%%%%%%%%%%%%%%%%%%%%%%%%%%%%%%%%%%%

\end{document}